\documentclass{siamltex}
\usepackage{lineno,hyperref}
\usepackage{multirow}
\usepackage[english]{babel}
\usepackage{amsmath}
\usepackage{bm}
\usepackage{ulem}
\usepackage[matha,mathx]{mathabx}
\usepackage{mathrsfs}
\usepackage{amssymb}
\usepackage{graphicx}
\usepackage{tikz}
\usepackage{here}
\usepackage{algorithm}
\usepackage{algorithmic}
\usepackage{graphicx}
\usepackage{amsmath}
\usepackage{rotating}
\usepackage{subfig}
\usepackage{here}
\usetikzlibrary{matrix,calc}
\usepackage{caption}
\usepackage{enumitem}

\makeatletter
\newcommand*{\rom}[1]{\expandafter\@slowromancap\romannumeral #1@}
\makeatother

\usepackage{tikz}
\usepackage{here}
\usepackage{algorithm}
\usepackage{algorithmic}
\usepackage{graphicx}   
\usepackage{amsmath}
\usepackage{here}
\usetikzlibrary{matrix,calc}
\usepackage{caption}
\DeclareMathAlphabet{\mathpzc}{OT1}{pzc}{m}{it}
\title{Non-negative Einstein tensor factorization for unmixing hyperspectral images}

\author{  A. El Hachimi\footnotemark[2] \thanks{The UM6P Vanguard Center, Mohammed VI Polytechnic University, Green 
City, Morocco.}
\and K. Jbilou\footnotemark[1] \thanks{Universit\'e du Littoral Cote d'Opale, LMPA, 50 rue F. Buisson, 62228 Calais-Cedex, France.}
\and A. Ratnani\footnotemark[1]}

\begin{document}
\maketitle 
\begin{abstract}
In this manuscript, we introduce a tensor-based approach to Non-Negative Tensor Factorization (NTF). The method entails tensor dimension reduction through the utilization of the Einstein product. To maintain the regularity and sparsity of the data, certain constraints are imposed. Additionally, we present an optimization algorithm in the form of a tensor multiplicative updates method, which relies on the Einstein product. To guarantee a minimum number of iterations for the convergence of the proposed algorithm, we employ the Reduced Rank Extrapolation (RRE) and the Topological Extrapolation Transformation Algorithm (TEA).  The efficacy of the proposed model is demonstrated through tests conducted on Hyperspectral Images (HI) for denoising, as well as for Hyperspectral Image Linear Unmixing. Numerical experiments are provided to substantiate the effectiveness of the proposed model for both synthetic and real data.
\end{abstract}

\section{Introduction}

Non-negative matrix factorization (NMF) serves as a method for dimensionality reduction; see the book of Nicolas Gillis \cite{gillis}. Given that real-world data is predominantly positive, the primary objective is to decompose a given non-negative matrix into the product of two non-negative matrices. The concept of NMF has drawn inspiration from Blind Signal Separation (BSS) \cite{ICA,BSS}. Let $A\in \mathbb{R}^{n\times m}$ be a non-negative matrix; i.e., $A_{i,j}\geq 0,\quad \forall (i,j)\in [n]\times [m]$, or it can be denoted also by $A\in \mathbb{R}^{n\times m}_+$, the NMF procedure involves determining the following decomposition
\[
A=XY + E,
\]
where $X\in \mathbb{R}^{n\times r}_+$, $Y\in \mathbb{R}^{r\times m}_+$, $E\in \mathbb{R}^{n\times m}$ represents the noise or the error, and $r$ is a positive integer such that $r\leq m$. $X$, $Y$, and $r$ may have different physical significations. In BSS problems, $X$ represents the mixing matrix, $Y$ plays the role of the source signal, and $r$ is the number of sources. In clustering problems, $X$ is considered as the basic matrix, while $Y$ denotes the weight matrix.

Images have emerged as one of the most significant technological advancements in recent years, finding applications across various domains such as medicine, geology, and archaeology. Initially, grayscale images were prevalent, followed by the invention of color images, which capture three bands (RGB): red, green, and blue. Additionally, 
Multispectral Imaging (MI) has been introduced as a method for evaluating images across a range of spectral bands, typically encompassing three or more colors. However, in many fields, including those requiring intricate analysis, color images, and MI may not suffice to achieve the desired precision. Consequently, the concept of Hyperspectral Images (HI) was introduced. Hyperspectral images break down an image into tens or hundreds of colors. With hundreds of contiguous spectral bands, hyperspectral images offer abundant spectral information, particularly beneficial for specific applications such as urban planning, object recognition, and mineral exploration. For more information, see \cite{book1,book_medecin,HI1}.

Among the interesting applications of NMF, we find Spectral mixture analysis, commonly referred to as unmixing, which stands as a widely utilized technique in remote sensing. This methodology is designed to disentangle mixed pixel spectra into constituent fractions associated with distinct pure spectral signatures, each accompanied by its corresponding proportion. Hyperspectral image unmixing (HU) encompasses a range of methods, including unsupervised, semi-supervised, and supervised approaches \cite{unsupervised1,unsupervised2,semisupervised,supervised}. Additionally, HU can be classified into geometric \cite{VCA} and statistical methods, with the NMF and NTF methods falling under the latter category \cite{Multi_updaes}.\\
Consider $a\in \mathbb{R}_+^{I}$ as the observed spectral pixel with $I$ wavelength-indexed bands. Linear mixing involves representing the $I$ wavelength-indexed bands using only $r$ end-member abundances. The formulation is provided by
\[
a=Xy + e,
\]
where the columns of the matrix $X\in \mathbb{R}_+^{I\times r}$ represent the endmembers ($r$ endmembers), $y\in \mathbb{R}_+^r$ is the vector of the abundance fractions for each end-member, and $e\in \mathbb{R}^I$ is an additive noise. While, when $J$ pixels are used, the formulation is transferred to matrix notation as
\[
A=X Y + E,
\]
where $X\in \mathbb{R}_+^{I\times r}$ is the matrix of endmember, $Y\in \mathbb{R}_+^{r\times J}$ is the abundance matrix, and $E\in \mathbb{R}^{I\times J}$ is an error. 
Additionally, two constraints must be imposed: the abundance non-negativity constraint (ANC), stipulating that the contribution from endmembers must be greater than or equal to zero, and the abundance sum-to-one constraint (ASC), requiring that the abundance fractions of all endmembers within a mixed pixel must sum up to one.\\
To approximate $X$ and $Y$, the following optimization problem can be considered
\[
\min_{X,Y} \dfrac{1}{2}\left\Vert A-XY\right\Vert_F^2,
\] 
and forcing the non-negativity of the obtained solutions. This problem may be solved by many approaches, but the most famous is the multiplicative updates approach given by Lee and Seung in \cite{Multi_updaes}.\\
Lately, there has been a growing focus on preserving sparsity. Given that sparsity is associated with the number of non-vanishing elements, the following optimization problem can be considered
\[
\min_{X,Y} \dfrac{1}{2}\left\Vert A-XY\right\Vert_F^2 + \lambda \left\Vert Y \right\Vert_0,\; \text{such that}\; X\geq 0, \, Y\geq 0,
\] 
with $\left\Vert Y\right\Vert_0=\texttt{card}\left\{ Y_{ij}/ Y_{ij}\neq 0,\; \text{for}\; i\in[r], \,j\in [J] \right\}$, in other words, it is the number of non-vanishing elements, it is called the $\ell_0$-norm (it is not a real norm), and $\lambda>0$ is a relaxation parameter.  While, this problem is known as an NP-hard problem \cite{l0-NP_hard1,l0-NP_hard2}. To turn around the difficulty of this optimization problem, the convexification of the $\ell_0$-norm is used, which is the $\ell_1$-norm that is given by 
\[
\left\Vert Y\right\Vert_1=\sum_{i=1}^{r}\sum_{j=1}^{J}\left\vert Y_{ij}\right\vert.
\] 
This gives the following optimization problem
\[
\min_{X,Y} \dfrac{1}{2}\left\Vert A-XY\right\Vert_F^2 + \lambda \left\Vert Y \right\Vert_1,\; \text{such that}\; X\geq 0, \, Y\geq 0.
\]
The aforementioned problem has been addressed using the multiplicative updates rule, as detailed in \cite{l1_NMF}. Additionally, numerous methods have been proposed to preserve sparsity, demonstrating enhanced efficiency. One such method involves the utilization of the $\ell_{1/2}$-norm, as discussed in \cite{l12_norm}.\\
Depending on the characteristics of the Hyperspectral image and the manner in which the data is mixed, regularization tools such as the $\ell_2$-norm may also be employed, as outlined in \cite{l2-NMF}.\\

Tensors have gained widespread adoption owing to their status as a generalization of scalars, vectors, and matrices. Notable among their attributes is the efficacity demonstrated in various applications. Tensors have found utility in a myriad of problem domains, including completion problems, denoising data, network analysis, and Hyperspectral image unmixing \cite{elha1,elha2,elhalouy,imbiriba2019low,mvnt}. Tensors can be approached through various logical frameworks, each of which is associated with a specific type of tensor product, such as the t-product, the c-product, and the n-mode product \cite{tproduct,cproduct,kolda}. In the context of tensors, several models for the NTF for HU have been introduced. Notable examples include the method based on block term decompositions (BTD) \cite{mvnt}, the work presented by Imbiriba et al. \cite{imbiriba2019low}, and the utilization of the n-mode product \cite{wu2023multidimensional}.\\
One of the notable tensor products is the Einstein product, which is widely recognized as one of the oldest \cite{einstein1,einstein2,einstein3}. It has been used frequently in many fields, for further examples, see \cite{tahiri,elha8}.

\medskip
The structure of this paper is outlined as follows: Section \ref{sec 1} provides an introduction to tensors and The Einstein product. Some theoretical results on Einstein NTF (ENTF) are illustrated in Section \ref{sec 33}. An approach for non-negative tensor decomposition based on the Einstein product is presented in Section \ref{sec 2}. Methods of extrapolation are employed in Section  \ref{sec extrapolation}, followed by the illustration of numerical experiments on Hyperspectral image denoising and Hyperspectral image linear unmixing in Section \ref{sec 3}.

\section{Preliminary}
\label{sec 1}
In this section, we delve into various concepts associated with tensors, beside an introduction onto the Einstein tensor product.

Over this paper, the notations used for tensors are the same as those employed by Kolda et al. \cite{kolda}. An $N$-order real tensor $\mathcal{A}$ is an $N$-array belongs to the space $\mathbb{R}^{I_1\times \ldots \times I_N}$. The $(i_1,\ldots,i_N)$-th element of $\mathcal{A}$ is denoted by $\mathcal{A}_{i_1\ldots i_N}$. In the sequel, the space of non-negative $N$-tensors of size $I_1\times \ldots \times I_N$ will be denoted by $\mathbb{R}^{I_1\times \ldots \times I_N}_+$.

\medskip

One of the renowned tensor products is the n-mode product, defined as follows.
\begin{definition}[\cite{kolda}]
	Let $\mathcal{A}\in\mathbb{R}^{I_1\times \ldots \times I_N}$ and $U\in \mathbb{R}^{J\times I_n}$ be an $N$-order tensor and a matrix. The resulting tensor $\mathcal{C}$ from the n-mode product between $\mathcal{A}$ and $U$ is given by 
	\[
	\mathcal{C}_{i_1\ldots i_{n-1}j i_{n+1}\ldots i_N}=\left(\mathcal{A}\times_n U\right)_{i_1\ldots i_{n-1} j  i_{n+1} \ldots i_N}=\sum_{i_n=1}^{I_n} \mathcal{A}_{i_1\ldots i_N}  U_{j i_n},
	\]
	where $\mathcal{C}$ is of size $\mathbb{R}^{I_1 \times \ldots \times I_{n-1}\times J\times I_{n+1}\times \ldots\times I_N}$.
\end{definition}

\medskip
\noindent A particular case of the n-mode product should be mentioned, which is when the product is between a tensor and a vector. For $\mathcal{A}\in \mathbb{R}^{I_1\times \ldots \times I_N}$ and $b\in \mathbb{R}^{I_n}$, then the n-mode product between $\mathcal{A}$ and $b$ is given by
\[
\left(\mathcal{A}\bar{\times}_n b\right)_{i_1 \ldots i_{n-1} i_{n+1} \ldots i_N}=\sum_{i_n=1}^{I_n} \mathcal{A}_{i_1\ldots i_N}  b_{i_n},
\]
and $\mathcal{A}\bar{\times}_n b\in \mathbb{R}^{I_1 \times \ldots \times I_{n-1}\times I_{n+1}\times \ldots\times I_N}.$

\medskip
The Einstein product \cite{einstein2} is given in the following definition.
\begin{definition}
	Let $\mathcal{A}\in \mathbb{R}^{I_1\times \ldots \times I_N \times J_1\times \ldots \times J_M}$ and $\mathcal{B}\in \mathbb{R}^{J_1\times \ldots \times J_M \times L_1 \times \ldots \times L_K}$, the Einstein product between $\mathcal{A}$ and $\mathcal{B}$ is the tensor $\mathcal{C}=\mathcal{A}*_M\mathcal{B}$  of size $I_1\times \ldots \times I_N \times L_1 \times \ldots \times L_K$, defined by 
	\[
	\left( \mathcal{A}*_M \mathcal{B}  \right)_{i_1 \ldots i_N \ell_1\ldots \ell_K}=\sum_{j_1,j_2,\ldots, j_M=1}^{J_1,J_2,\ldots,J_M} \mathcal{A}_{i_1 \ldots i_N j_1\ldots j_M} \mathcal{B}_{j_1 \ldots j_M \ell_1\ldots \ell_K}.
	\]
\end{definition}

\medskip
\noindent If $\mathcal{A}\in \mathbb{R}^{I_1\times \ldots \times I_N\times J}$ be an $(N+1)$-th order tensor and $\mathcal{B}\in \mathbb{R}^{J\times I}$, then the n-mode product of $\mathcal{A}$ and $B^T$ over the mode $N+1$ is the same as the Einstein product when summing over the last mode of the tensor, i.e., 
\[
\mathcal{A}*_1B=\mathcal{A}\times_{N+1} B^T,
\]
also, if $\mathcal{A}\in \mathbb{R}^{J\times I_1 \times \ldots \times \dots I_N}$ and $B^{I\times J}$, then we can write 
\[
B *_1 \mathcal{A}= \mathcal{A}\times_1 B^T.
\]
\medskip

Several particular tensors associated with the Einstein product \cite{einstein2, einstein3} are given as follows:
\begin{itemize}
	\item{\bf Transpose tensor:}  For a given tensor $\mathcal{A}\in \mathbb{R}^{I_1\times \ldots \times I_N \times J_1 \times \ldots \times J_M}$, the tensor $\mathcal{B}\in \mathbb{R}^{J_1\times \ldots \times J_M \times I_1 \times \ldots \times I_N}$ is defined as the transpose of $\mathcal{A}$ if $\mathcal{B}_{j_1\ldots j_M i_1\ldots i_N}=\mathcal{A}_{i_1\ldots i_N j_1 \ldots j_M}$. This transpose tensor is denoted by $\mathcal{A}^T$.
	\item{\bf Diagonal tensor:} A diagonal tensor $\mathcal{D}\in \mathbb{R}^{I_1\times \ldots \times I_N \times J_1\times \ldots \times J_M}$ is a tensor that verifies $\mathcal{D}_{i_1\ldots i_N j_1\dots j_M}=0$ if $i_k\neq j_k$ for all $k\in \{1,\ldots, \min(N,M)\}$.
	\item{\bf Identity tensor:} $\mathcal{I}_N\in \mathbb{R}^{I_1\times \ldots \times I_N \times I_1\times \ldots \times I_N}$ is said to be an identity tensor under the Einstein product if it is diagonal and the non-vanishing elements are equal to $1$.
	\item{\bf Orthonormal tensor:} $\mathcal{Q}\in \mathbb{R}^{I_1\times \ldots \times I_N \times I_1 \times \ldots \times I_N}$ is said to be an orthogonal tensor if 
	\[
	\mathcal{Q}^T*_N \mathcal{Q}=\mathcal{Q}*_N \mathcal{Q}^T=\mathcal{I}_N.
	\]
	\item{\bf  Tensor inverse:} The tensor $\mathcal{A}\in \mathbb{R}^{I_1\times \ldots \times I_N\times I_1\times \ldots \times I_N}$ is invertible if there exists a tensor $\mathcal{B}\in \mathbb{R}^{I_1\times \ldots \times I_N \times I_1 \times \ldots \times I_N}$ such that 
	\[
	\mathcal{A}*_N \mathcal{B}=\mathcal{B}*_N \mathcal{A}=\mathcal{I}_N,
	\]  
	if this tensor exists, it is denoted $\mathcal{A}^{-1}$.
\end{itemize}

\medskip
Consider two tensors $\mathcal{A}$ and $\mathcal{B}$ of the same size $I_1 \times \ldots \times I_N \times J_1 \times \ldots \times J_M$. The inner product between $\mathcal{A}$ and $\mathcal{B}$ is given by
\[
\langle \mathcal{A}, \mathcal{B} \rangle={\bf tr}\left( \mathcal{A}^T *_N \mathcal{B}  \right)={\bf tr}\left( \mathcal{B}^T *_N \mathcal{A}  \right),
\]
with
\[
{\bf tr}\left(\mathcal{A}\right)=\sum_{i_1,\ldots i_N=1}^{I_1,\ldots, I_N} \mathcal{A}_{i_1\ldots i_N i_1 \ldots i_N}.
\]
The corresponding norm is given by 
\[
\left\Vert \mathcal{A}\right\Vert_F=\sqrt{{\bf tr}\left(\mathcal{A}^T *_N \mathcal{A}\right)}.
\]

\medskip
\begin{proposition}
	Let $\mathcal{A}\in \mathbb{R}^{I_1\times \ldots \times I_N\times J_1 \times \ldots \times J_M}$ and $\mathcal{B}\in \mathbb{R}^{J_1\times \ldots \times J_M \times L_1 \times \ldots \times L_K}$, then 
	\[
	\left(\mathcal{A}*_M \mathcal{B}\right)^T = \mathcal{B}^T *_M \mathcal{A}^T.
	\]
\end{proposition}
\medskip
\begin{proposition}
	Let $\mathcal{A}\in \mathbb{R}^{I_1 \times \ldots \times I_N \times J_1 \times \ldots \times J_M}$ and $\mathcal{U}\in \mathbb{R}^{J_1\times \ldots \times J_M\times J_1\times \ldots \times J_M}$, where $\mathcal{U}$ is an orthonormal tensor. Then 
	\[
	\left\Vert \mathcal{A}*_M \mathcal{U}\right\Vert_F= \left\Vert \mathcal{A}\right\Vert_F.
	\]
\end{proposition}

\noindent The extension of the Singular Value Decomposition (SVD) to tensors using the Einstein product is provided below \cite{einstein2}.
\begin{theorem}
	The tensor singular value decomposition of $\mathcal{A}\in \mathbb{R}^{I_1\times \ldots \times I_N\times J_1 \times \ldots \times J_M}$.  $\mathcal{A}$ is given by 
	\[
	\mathcal{A}=\mathcal{U}*_N \mathcal{S}*_M \mathcal{V}^T,
	\]
	where $\mathcal{U}\in \mathbb{R}^{I_1\times \ldots \times I_N\times I_1\times \ldots \times I_N}$ and $\mathcal{V}\in \mathbb{R}^{J_1 \times \ldots \times J_M\times J_1 \times \ldots \times J_M}$ are orthonormal tensors and $\mathcal{S}\in \mathbb{R}^{I_1\times \ldots \times I_N\times J_1 \times \ldots \times J_M}$ is a diagonal tensor.
\end{theorem}

\medskip
\noindent Assume that $M\leq N$. As has been discussed in \cite{elha8}, the Tensor SVD for a tensor $\mathcal{A}\in \mathbb{R}^{I_1\times \ldots \times I_N \times J_1 \times \ldots \times J_M}$ can be also written as 
\[
\mathcal{A}=\sum_{i_1,\ldots, i_M=1}^{\widetilde{I}_1, \ldots, \widetilde{I}_M}\mu_{i_1\ldots i_M}\mathcal{U}_{i_1\ldots i_M}\circ \left(   \mathcal{V}_{i_1\ldots i_M}  \right)^T,
\]
where $\widetilde{I}_i=\min\left(I_i,J_i\right)$ for $i=1,\ldots M$, and $\mu_{i_1\ldots i_M}$, $\mathcal{U}_{i_1\ldots i_M}$, and $\mathcal{V}_{i_1\ldots i_M}$ are respectively the singular value, left and right singular tensor associated to $\mathcal{A}$, and where $\circ$ represents the outer product. Where $\mu_{i_1\ldots i_M}$'s are ordered decreasingly, i.e., 
\[
\mu_{1\ldots 1}\geq \mu_{1\, 2\ldots 1} \geq \ldots \geq \mu_{\widetilde{I}_1\ldots \widetilde{I}_M}\geq 0.
\]

\noindent The number $\texttt{card}\left\{  \mu_{i_1\ldots i_M}\neq 0 / i_j\in [I_j], \; \text{for } j=1,\ldots,M \right\}$ is called the rank of $\mathcal{A}$.

\medskip
\begin{definition}
	Let $\mathcal{A}\in \mathbb{R}^{I_1 \times \ldots \times I_N \times J_1 \times \ldots \times J_M}$.
	\begin{itemize}
		\item If ${\tt rank}\left( \mathcal{A} \right)=J_1\ldots J_M$ with $J_1\ldots J_M\leq I_1 \ldots I_N$. Then $\mathcal{A}$ has a left inverse denoted by $\mathcal{A}_{left}^{-1}$ of size $J_1\times \ldots \times J_M\times I_1 \times \ldots \times I_N$ such that 
		\[
		\mathcal{A}_{left}^{-1} *_N\mathcal{A}=\mathcal{I}_M.
		\]
		\item If ${\tt rank}\left( \mathcal{A} \right)=I_1\ldots I_M$ with $I_1\ldots I_N\leq J_1 \ldots J_M$. Then $\mathcal{A}$ has a right inverse denoted by $\mathcal{A}_{right}^{-1}$ of size $J_1\times \ldots \times J_M\times I_1 \times \ldots \times I_N$ such that 
		\[
		\mathcal{A}*_M\mathcal{A}_{right}^{-1}=\mathcal{I}_N.
		\]
	\end{itemize}
\end{definition}

A simple generalization of the Gram-Schmidt process (QR decomposition) is going to be proposed to be used later. This process consists in orthonormalizing a sequence of tensors $\mathcal{A}_1,\ldots, \mathcal{A}_k$ in $\mathbb{R}^{I_1\times \ldots \times I_N}$. If we consider the tensor $\mathcal{A}=\left[ \mathcal{A}_1,\ldots, \mathcal{A}_k \right]\in \mathbb{R}^{I_1\times \ldots \times I_N\times k}$, we are going to find the following decomposition
\[
\mathcal{A}=\mathcal{Q}\times_{N+1} R^T,
\]
with $\mathcal{Q}\in \mathbb{R}^{I_1\times \ldots \times I_N\times k}$ and $R\in \mathbb{R}^{k\times k}$, such that if we write $\mathcal{Q}=\left[\mathcal{Q}_1,\ldots, \mathcal{Q}_k  \right]$, then we find $\langle \mathcal{Q}_i, \mathcal{Q}_j \rangle=\delta_{i,j}$, or otherwise $\mathcal{Q}^T *_N\mathcal{Q}=I_{k}$, and $R$ is an upper triangular matrix.

\begin{algorithm}[H]
	\caption{Tensor Einstein QR decomposition (E-QR)}
	\label{alg E-QR}
	\textbf{Input:} $\mathcal{A}=\left[\mathcal{A}_1, \ldots, \mathcal{A}_k\right]\in \mathbb{R}^{I_1\times \ldots \times I_N\times k}$.\\
	\textbf{Output:} $\mathcal{Q}=\left[\mathcal{Q}_1, \ldots, \mathcal{Q}_k\right]\in \mathbb{R}^{I_1\times \ldots \times I_N\times k}$ and $R$, such that $\mathcal{A}=\mathcal{Q}\times_{N+1} R^T$.
	\begin{algorithmic}[1]
		\STATE Compute $R_{1,1}=\left\Vert \mathcal{A}_1\right\Vert_F$, $\mathcal{Q}_1=\dfrac{\mathcal{A}_1}{R_{1,1}}$.\\
		\FOR{$j=2:k$}
		\STATE $\mathcal{W}=\mathcal{A}_j$
		\FOR{$i=1:j-1$}
		\STATE $R_{i,j}=\langle \mathcal{Q}_i, \mathcal{W} \rangle $; 
		\STATE $\mathcal{W}=\mathcal{W}-R_{i,j}\mathcal{Q}_i$.
		\ENDFOR
		\STATE $R_{j,j}=\left\Vert \mathcal{W}\right\Vert_F$.
		\STATE $\mathcal{Q}_j=\dfrac{\mathcal{W}}{R_{j,j}}$.
		\ENDFOR
	\end{algorithmic}
\end{algorithm}
\medskip

\section{NTF based on the Einstein product (ENTF)}
\label{sec 33}
In this section, we define the NTF based on the Einstein product (ENTF), besides a geometrical interpretation, and some related properties.
\medskip
\begin{definition}[ENTF]
	Let $\mathcal{A}\in \mathbb{R}_+^{I_1\times \ldots \times I_N\times J_1 \times \ldots \times J_M}$ be an $(N+M)$-order non-negative tensor. The $r\in \mathbb{N}^*$ ENFT consists in computing if possible two tensors $\mathcal{X}\in \mathbb{R}_+^{I_1\times \ldots \times I_N\times r}$ and $\mathcal{Y}\in \mathbb{R}_+^{r\times J_1 \times \ldots \times J_M}$ such that
	\begin{equation}
		\mathcal{A}=\mathcal{X}*_1 \mathcal{Y}.\label{eq NTFE}
	\end{equation}
\end{definition}

\medskip
\noindent Equation \eqref{eq NTFE} can be formulated by using the outer product, i.e., 
\[
\mathcal{A}=\sum_{i=1}^{r} \mathcal{X}(:,\ldots, :, i)\circ \mathcal{Y}(i,:,\ldots,:),
\] 
The above decomposition can also be called an ENTF of size $r$.
\medskip
\begin{definition}[Einstein non-negative rank]
	The Einstein non-negative rank of a tensor $\mathcal{A}\in \mathbb{R}_+^{I_1\times \ldots \times I_N \times J_1 \times \ldots \times J_M}$ is the smallest number $r$ such that the factorization in \eqref{eq NTFE} admits an exact $r$ ENTF. It is denoted by $\texttt{rank}_+\left(\mathcal{A}\right)$.
\end{definition}

\medskip
\noindent Obviously, it is easy to see that all time ${\tt rank}\left(\mathcal{A}\right)\leq {\tt rank}_+\left(\mathcal{A}\right)\leq \min(I_1\ldots I_N, J_1\ldots J_M)$.

\medskip
\begin{definition}
	Let $\mathcal{A}\in \mathbb{R}^{J_1 \times \ldots  \times J_M \times I_1 \times \ldots \times I_N }$. We define the following spaces.
	\begin{itemize}
		\item The ${\tt cone}\left(\mathcal{A}\right)$ is the cone  of $\mathcal{A}$.
		\begin{align*}
			{\tt cone}\left(\mathcal{A}\right)=\left\{ \mathcal{X}\in \mathbb{R}^{J_1\times \ldots \times J_M}/  \mathcal{X}=\mathcal{A}*_N \mathcal{Y}, \quad \mathcal{Y}\in \mathbb{R}^{I_1\times \ldots \times I_N}_+ \right\}.
		\end{align*}
		\item The ${\tt Im}\left(\mathcal{A}\right)$ is the image of $\mathcal{A}$.
		\begin{align*}
			{\tt Im}\left(\mathcal{A}\right)=\left\{ \mathcal{X}\in \mathbb{R}^{J_1\times \ldots \times J_M}/   \mathcal{X}=\mathcal{A}*_N \mathcal{Y}, \quad \mathcal{Y}\in \mathbb{R}^{I_1\times \ldots \times I_N} \right\}.
		\end{align*}
		\item The convex hull of $\mathcal{A}$ is defined as
		\begin{align*}
			\texttt{conv}\left( \mathcal{A} \right)=&\left\{ \mathcal{X}/ \mathcal{X}=\mathcal{A}*_N \mathcal{Y}, \quad \mathcal{Y}\in \mathbb{R}_+^{I_1 \times \ldots \times I_N},\;  \textbf{ones}_{I_1,\ldots , I_N}*_N\mathcal{Y}=1  \right\},
		\end{align*}
		where $\textbf{ones}_{I_1,\ldots , I_N}$ is the tensor of size $I_1\times \ldots \times I_M$ such that all its elements are ones.
		\item The unit simplex $\Delta_{I_1, \ldots , I_N}$
		\[
		\Delta_{I_1, \ldots , I_N}=\left\{ \mathcal{X}\in \mathbb{R}_+^{I_1\times \ldots \times J_N}/ \;  \textbf{ones}_{I_1,\ldots, I_N}*_N \mathcal{X}=1  \right\}.
		\]
		\item The affine hull of $\mathcal{A}$
		\[
			{\tt aff}\left(\mathcal{A}\right)=\left\{  \mathcal{X}/\, \mathcal{X}=\mathcal{A}*_N \mathcal{Y}, 
 \, \mathcal{Y}\in \mathbb{R}^{I_1\times \ldots \times I_N},\, \textbf{ones}_{I_1,\ldots, I_N}*_N \mathcal{Y}=1  \right\}.
		\]
	\end{itemize}
\end{definition}

\medskip
Assume we have an Exact ENTF for $\mathcal{A}\in \mathbb{R}^{I_1 \times \ldots \times I_N \times J_1 \times \ldots \times J_M}$ of size $r$; $\mathcal{A}=\mathcal{X}*_1 \mathcal{Y}$. Since 
\[
\mathcal{A}(:,\ldots, :, j_1, \ldots, j_M)=\mathcal{X}*_1 \mathcal{Y}(:,j_1, \ldots, j_M),
\]
then, it can be seen that 
\[
\mathcal{A}(:,\ldots, :, j_1, \ldots, j_M)\in {\tt cone}\left(\mathcal{X}\right)\subset \mathbb{R}^{I_1\times \ldots \times I_N}.
\]
This is equivalent to 
\[
{\tt cone}\left( \mathcal{A}  \right)\subset {\tt cone}\left(\mathcal{X}\right) \subset \mathbb{R}^{I_1\times \ldots \times I_N}.
\]
Therefore, for a given non-negative tensor $\mathcal{A}\in \mathbb{R}^{I_1 \times \ldots \times I_N\times J_1 \times \ldots \times J_M}_+$, finding an ENTF is equivalent to find a tensor $\mathcal{X}$, such that its convex cone is between the two spaces ${\tt cone}\left(\mathcal{A}\right)$ and $\mathbb{R}^{I_1\times \ldots \times I_N}$. Assume we have the following ENTF, 
\begin{equation}
	\label{eq 3.3}
	\begin{split}
		\mathcal{A}=\mathcal{X}*_1 \mathcal{Y} \Longleftrightarrow & \mathcal{A}*_M\mathcal{D}_\mathcal{A}=\left(\mathcal{X}*_1\mathcal{D}_\mathcal{X}\right)*_1 \left( \mathcal{D}_\mathcal{X}^{-1}*_1\mathcal{Y}*_M\mathcal{D}_\mathcal{A} \right),\\
		\Longleftrightarrow & \theta\left(  \mathcal{A} \right)= \theta \left( \mathcal{X}  \right)*_1 \mathcal{Y}',
	\end{split}
\end{equation}
where $\mathcal{Y}'=\mathcal{D}_\mathcal{X}^{-1}*_1\mathcal{Y}*_M\mathcal{D}_\mathcal{A}$, $\theta\left(\mathcal{A}\right)=\mathcal{A}*_M \mathcal{D}_\mathcal{A}$, and $\theta\left(\mathcal{X}\right)=\mathcal{X}*_1 \mathcal{D}_\mathcal{X}$, with
\[
\left(\mathcal{D}_\mathcal{A}\right)_{i_1 \ldots i_N j_1 \ldots j_M}=\begin{cases} \dfrac{1}{\left\Vert \mathcal{A}_{:\ldots : j_1\ldots j_M}\right\Vert_1},& \text{if}\quad i_k=j_k;\, k=1,\dots,K,\\
	0 & \text{otherwise}.
\end{cases}
\]
with $K=\min(N,M)$, and
\[
\left(\mathcal{D}_\mathcal{X}\right)_{i_1 \ldots i_N j}=\begin{cases} \dfrac{1}{\left\Vert \mathcal{X}_{:, 1\ldots 1,j}\right\Vert_1},& \quad \text{if}\quad i_1=j,\\
	0 & \text{otherwise}.
\end{cases}
\]
\noindent It can be seen that $\theta\left( \mathcal{A}  \right)\in {\tt conv}\left( \theta\left(\mathcal{X}\right)  \right)$. From \eqref{eq 3.3} we can conclude
\[
\texttt{cone}\left(\mathcal{A}\right)\subset \texttt{cone}\left(\mathcal{X}\right)\subset \mathbb{R}^{I_1\times \ldots \times I_N}, \]
is equivalent to
\[ \texttt{conv}\left(\theta\left(\mathcal{A}\right)\right)\subset \texttt{conv}\left(\theta\left(\mathcal{X}\right)\right)\subset \mathbb{R}^{I_1\times \ldots \times I_N}.
\]

\medskip
In the sequel of this section, we will provide some cases of the non-negative rank of a tensor.
\begin{lemma}
	Let $\mathcal{A}\in \mathbb{R}^{I_1\times \ldots \times I_N\times J_1\times \ldots \times J_M}_+$ be an $(N+M)$-th order tensor. Then,  $\texttt{conv}\left(\theta\left(\mathcal{A}\right)\right)$ and $\texttt{Im}\left( \mathcal{A} \right)\bigcap \Delta_{I_1,\ldots, I_N}$ are of dimension $\texttt{rank}\left(\mathcal{A}\right)-1$. 
\end{lemma}
\medskip

\begin{theorem}
	\label{theo 3.5}
	Assume that $\mathcal{A}\in \mathbb{R}^{I_1\times \ldots \times I_N \times J_1\times \ldots \times J_M}_+$. Assume that ${\tt rank}\left(\mathcal{A}\right)=2$, then ${\tt rank}\left(\mathcal{A}\right)={\tt rank}\left(\mathcal{A}\right)_+$.
\end{theorem}

\medskip
\begin{proof}
	The proof of Lemma and Theorem are straightforward using the matrix case.
\end{proof}

\medskip
\begin{lemma}
	Let $\mathcal{A}\in \mathbb{R}^{I_1\times \ldots \times I_N \times J_1\times \ldots \times J_M}_+$ and $\mathcal{A}=\mathcal{X}*_1 \mathcal{Y}$ be an exact ENTF of size $r={\tt rank}\left(\mathcal{A}\right)$. Then ${\tt rank}\left(\mathcal{X}\right)={\tt rank}\left(\mathcal{A}\right)$ and $Im\left(\mathcal{X}\right)=Im\left(\mathcal{A}\right)$.
\end{lemma}

\medskip

\begin{proof}
	Since $\mathcal{A}=\mathcal{X}*_1 \mathcal{Y}$ is an exact ENTF of size $r$. Then 
	\begin{eqnarray*}
		r={\tt rank}\left(\mathcal{A}\right)={\tt rank}\left( \mathcal{X}*_1 \mathcal{Y} \right)&\leq& \min\left(  {\tt rank}\left(\mathcal{X}\right), {\tt rank}\left(\mathcal{Y}\right)\right)\\
		&\leq& {\tt rank}\left(\mathcal{X}\right)\\
		&\leq & \min\left( I_1\ldots I_N, r   \right)\\
		&\leq& r.
	\end{eqnarray*}
	Consequently $r={\tt rank}\left(\mathcal{A}\right)={\tt rank}\left(\mathcal{X}\right)$.\\
	The last equality is trivial.
\end{proof}

\medskip

\noindent The following definition discusses when an ENTF is said to be unique.
\begin{definition}[Uniqueness of ENTF]
	The exact ENTF $\mathcal{A}=\mathcal{X}*_1 \mathcal{Y}$ of size $r$ is unique if and only if for any other exact ENTF $\mathcal{A}=\mathcal{X}'*_1 \mathcal{Y}'$ of size $r$, there exists a permutation matrix $P\in \mathbb{R}^{r\times r}$ and a diagonal matrix $D$ with positive diagonal elements such that 
	\[
	\mathcal{X}' = \mathcal{X} \times_{N+1} \left( D P^T\right) \quad \text{and}\quad 	\mathcal{Y}'=\mathcal{Y} \times_{1} \left( P D^{-1}\right).
	\]
\end{definition}

\medskip
\noindent The following theorem provides some properties of the non-negative rank of a tensor.
\begin{theorem}
	Let $\mathcal{A}\in \mathbb{R}^{I_1\times \ldots \times I_N\times J_1\times \ldots \times J_M}_+$. The following properties are satisfied.
	\begin{enumerate}
		\item $\texttt{rank}_+\left(\mathcal{A}\right)=\texttt{rank}_+\left(\mathcal{A}^T\right)$.
		\item For all $\mathcal{B}$ in $\mathbb{R}^{I_1\times \ldots \times I_N\times J_1\times \ldots \times J_M}_+$, $\texttt{rank}_+\left(\mathcal{A}+\mathcal{B}\right)\leq \texttt{rank}_+\left(\mathcal{A}\right)+\texttt{rank}_+\left(\mathcal{B}\right)$.
		\item If $\texttt{rank}_+\left(\mathcal{A}\right)\leq 2$ or $\min(I_1\ldots I_N, J_1\ldots J_M)={\tt rank}\left(\mathcal{A}\right)$, then
		\[
		\texttt{rank}_+\left(\mathcal{A}\right)=	\texttt{rank}\left(\mathcal{A}\right).
		\]
	\end{enumerate}
\end{theorem}

\medskip
\begin{proof}
	\begin{enumerate}
		\item Let $\mathcal{A}=\mathcal{X}*_1 \mathcal{Y}$ be an exact ENTF. Then $\mathcal{A}^T=\mathcal{Y}^T*_1 \mathcal{X}^T$. Then $\texttt{rank}_+\left(\mathcal{A}\right)=\texttt{rank}_+\left(\mathcal{A}^T\right)$.
		\item Let $\mathcal{A}=\mathcal{X}*_1\mathcal{Y}$ by an exact ENTF of rank $r_\mathcal{A}$, and $\mathcal{B}=\mathcal{Z}*_1\mathcal{W}$ by an exact ENTF of rank $r_\mathcal{B}$. \\
		We can see that 
		\[
		\mathcal{A}+\mathcal{B}=\mathcal{X}*_1\mathcal{Y}+\mathcal{Z}*_1\mathcal{W}=\begin{pmatrix}
			\mathcal{X} \; \mathcal{Z}
		\end{pmatrix}*_1\begin{pmatrix}
			\mathcal{Y}\\
			\mathcal{W}
		\end{pmatrix},
		\]
		which provides a sum of $r_\mathcal{A}+r_\mathcal{B}$ elements. Which means that as least, the rank of $\mathcal{A}+\mathcal{B}$ is $r_\mathcal{A}+r_\mathcal{B}$, in other word
		\[
		\texttt{rank}_+\left(\mathcal{A}+\mathcal{B}\right)\leq \texttt{rank}_+\left(\mathcal{A}\right)+\texttt{rank}_+\left(\mathcal{B}\right).
		\]
		\item If 	$\texttt{rank}_+\left(\mathcal{A}\right)\leq 2$, it has been proved in Theorem \ref{theo 3.5}.\\
		If $\min(I_1\ldots I_N, J_1\ldots J_M)={\tt rank}\left(\mathcal{X}\right)$, then we will find $\mathcal{A}=I_{I_1}*_1 \mathcal{X}=\mathcal{X}*_1 I_{J_M}$. Thus $\texttt{rank}_+\left(\mathcal{A}\right)=	\texttt{rank}\left(\mathcal{A}\right)$. ($I_{I_1}$ is the matrix identity of size $I_1\times I_1$, and similarly for $I_{J_M}$).
	\end{enumerate}
\end{proof}

\medskip
\begin{lemma}
	Let $\mathcal{X}\in \mathbb{R}_+^{I_1\times \ldots \times I_N \times r}$, and $\mathcal{Y}\in \mathbb{R}_+^{r\times J_1\times \ldots \times J_M}$, be the factors of an exact $r$ ENTF of $\mathcal{A}\in \mathbb{R}_+^{I_1\times \ldots \times I_N \times J_1 \times \ldots \times J_M}$, with $r={\tt rank}\left(\mathcal{X}\right)$. Then, any factors $\mathcal{X}'$ and $\mathcal{Y}'$ of an exact $r$ ENTF of $\mathcal{X}$ verifies
	\[
	\mathcal{X}'= \mathcal{X}\times_{N+1} Q\geq 0 \quad \text{and}\quad \mathcal{Y}'= \mathcal{Y}\times_1 Q^{-1} \geq 0,
	\]
	for some invertible matrix $Q$ of size $r\times r$.
\end{lemma}

\medskip
\begin{proof}
	Since $\mathcal{A}=\mathcal{X}*_1 \mathcal{Y}=\mathcal{X}'*_1\mathcal{Y}'$ are exact ENTF with size $r=\texttt{rank}\left(\mathcal{A}\right)$. Then we have
	\[
	Im\left(\mathcal{A}\right)=Im\left(\mathcal{X}\right)=Im\left(\mathcal{X}'\right),
	\]
	Then there exists an invertible matrix $Q$ of size $r\times r$ such that 
	\[
	\mathcal{X}'=\mathcal{X}\times_{N+1} Q.
	\]
	Moreover, since $\mathcal{X}$ is of rank $r$. Then, it admits a left inverse $\mathcal{X}^{-1}_{\text{left}}$ such that, $\mathcal{X}^{-1}_{\text{left}}*_N \mathcal{X}=I_r$. Thus
	\[
	\mathcal{X}'*_1 \mathcal{Y}'=\mathcal{X}\times_{N+1} Q *_1 \mathcal{Y}'=\mathcal{X}*_1 \mathcal{Y}.
	\]
	Consequently
	\[
	\mathcal{Y}'\times_1 Q=\mathcal{Y} \Leftrightarrow  \mathcal{Y}'=\mathcal{Y}\times_1 Q^{-1}.
	\]
\end{proof}
\medskip
\begin{theorem}
	The Exact $r$ ENTF of factors $\mathcal{X}$ and $\mathcal{Y}$ of a tensor $\mathcal{X}\in \mathbb{R}_+^{I_1\times \ldots \times I_N \times J_1 \times \ldots \times J_M}$, with $r={\tt rank}\left(\mathcal{X}\right)$ is unique if and only if the only invertible matrices $\mathcal{Q}$ such that
	\[
	\mathcal{X}\times_{N+1}  Q\geq 0, \quad \text{and}\quad  \mathcal{Y}\times_1 Q^{-1}\geq 0,
	\]
	are diagonal matrices with positive diagonal elements, up to perturbations of rows and columns.
\end{theorem}
\medskip
\begin{proof}
	The proof follows from the last Lemma and the fact that both $Q$ and $Q^{-1}$ are non-negative if and only if $Q$ is diagonal with positive diagonal elements, up to perturbations of rows and columns.
\end{proof}
\section{NTF by the Einstein product based on HI}
\label{sec 2}
Let $\mathcal{A}\in \mathbb{R}^{I\times J \times K}$ be a hyperspectral image, where the first index represents the number of bands ($I$ bands), while the second and the third represent the spatial index. The NTF-method based on the Einstein product consists of factorizing $\mathcal{A}$ into 
\[
\mathcal{A}=\mathcal{X}*_1 \mathcal{Y},
\]
such that $\mathcal{X}\in \mathbb{R}^{I\times r}_+$, $\mathcal{Y}\in \mathbb{R}^{r\times J\times K}_+$.\\
To this end, the following optimization problem can aid in approximating this factorization
\begin{equation}
	\begin{split}
		&\min_{\mathcal{X}, \mathcal{Y}} \dfrac{1}{2}\left\Vert \mathcal{A} - \mathcal{X}*_1 \mathcal{Y}\right\Vert_F^2,\\
		&\text{subject that } \mathcal{X}\geq 0,\; \mathcal{Y}\geq 0,\;  \mathcal{Y}\bar{\times}_1 {\bf 1}_r={\bf ones}_{J,K},
	\end{split}
\end{equation}
where $\bf{1}_r$ is the vector of size $r$ all its elements are ones; $\bf{1}_r=ones(r,1)$, and ${\bf ones}_{J,K}\in \mathbb{R}^{J\times K}$ is a matrix that all its elements are ones.  It is too hard to keep the sparsity constraint for sparse data, for this reason, some tools are used. Among them, there is the $\ell_0$-norm. So the optimization problem will be as follows
\[
\begin{split}
	&\min_{\mathcal{X}, \mathcal{Y}} \dfrac{1}{2}\left\Vert \mathcal{A} - \mathcal{X}*_1 \mathcal{Y}\right\Vert_F^2 + \lambda_s \left\Vert \mathcal{Y}\right\Vert_0,\\
	&\text{subject that } \mathcal{X}\geq 0,\; \mathcal{Y}\geq 0,\; \mathcal{Y}\bar{\times}_1 {\bf 1}_r={\bf ones}_{J,K},
\end{split}
\]
where $\lambda_s>0$ is a parameter linked with the sparsity of the data.\\
The problem with using the $l_0$-norm is that it is an NP-hard problem. To turn around this problem, the relaxation of the $\ell_0$-norm is used. The relaxation of the $\ell_0$-norm is the $\ell_1$-norm, thus the approximation of the desired solution could be obtained by solving the following optimization problem
\[
\begin{split}
	&\min_{\mathcal{X}, \mathcal{Y}} \dfrac{1}{2}\left\Vert \mathcal{A} - \mathcal{X}*_1 \mathcal{Y}\right\Vert_F^2 + \lambda_s \left\Vert \mathcal{Y}\right\Vert_1,\\
	&\text{subject that } \mathcal{X}\geq 0,\; \mathcal{Y}\geq 0,\; \mathcal{Y}\bar{\times}_1 {\bf 1}_r={\bf ones}_{J,K}.
\end{split}
\]
For some regularization reasons, we can reformulate the above optimization problem to be as
\begin{equation*}
	\begin{split}
		&\min_{\mathcal{X}, \mathcal{Y}, \mathcal{P}, \mathcal{Q}} \dfrac{1}{2}\left\Vert \mathcal{A} - \mathcal{X}*_1 \mathcal{Y}\right\Vert_F^2 + \lambda_s\left\Vert \mathcal{Y}\right\Vert_1 + \dfrac{\lambda_\mathcal{Y}}{2}\left\Vert \mathcal{Y} - \mathcal{P}\right\Vert_F + \dfrac{\lambda_\mathcal{X}}{2}\left\Vert \mathcal{X}-\mathcal{Q}\right\Vert_F,\\
		&\text{subject that } \mathcal{X}\geq 0,\; \mathcal{Y}\geq 0,\; \mathcal{Y}\bar{\times}_1 {\bf 1}_r={\bf ones}_{J,K},
	\end{split}
\end{equation*}
where $\mathcal{P}\in \mathbb{R}^{r\times J\times K}$ and $\mathcal{Q}\in \mathbb{R}^{I\times r}$ are quantities approximating $\mathcal{Y}$ and $\mathcal{X}$, respectively, and $\lambda_\mathcal{Y}>0$ and $\lambda_\mathcal{X}>0$ are relaxation parameters.\\
Since we aim to find a solution $\mathcal{Y}$ close to $\mathcal{P}$,then we can rewrite our optimization problem as
\[
\begin{split}
	&\min_{\mathcal{X}, \mathcal{Y}, \mathcal{P}, \mathcal{Q}} \dfrac{1}{2}\left\Vert \mathcal{A} - \mathcal{X}*_1 \mathcal{Y}\right\Vert_F^2 + \lambda_s\left\Vert \mathcal{P}\right\Vert_1 + \dfrac{\lambda_\mathcal{Y}}{2}\left\Vert \mathcal{Y} - \mathcal{P}\right\Vert_F+ \dfrac{\lambda_\mathcal{X}}{2}\left\Vert \mathcal{X}-\mathcal{Q}\right\Vert_F,\\
	&\text{subject that } \mathcal{X}\geq 0,\; \mathcal{Y}\geq 0,\; \mathcal{Y}\bar{\times}_1 {\bf 1}_r={\bf ones}_{J,K}.
\end{split}
\] 
Then, we get 
\begin{eqnarray}
	\mathcal{X}_{k+1}&=&\arg\min_{\mathcal{X}} \dfrac{1}{2} \left\Vert \mathcal{A}- \mathcal{X}*_1 \mathcal{Y}_k \right\Vert_F^2 + \dfrac{\lambda_\mathcal{X}}{2}\left\Vert \mathcal{X}-\mathcal{Q}\right\Vert_F.\label{eq X}\\
	\mathcal{Y}_{k+1}&=&\arg\min_{\mathcal{Y}} \dfrac{1}{2} \left\Vert \mathcal{A}- \mathcal{X}_{k+1}*_1 \mathcal{Y} \right\Vert_F^2  + \dfrac{\lambda_\mathcal{Y}}{2}\left\Vert \mathcal{Y}-\mathcal{P}\right\Vert_F.\label{eq Y}\\
	\mathcal{P}_{k+1}&=&\arg\min_{\mathcal{P}}  \dfrac{\lambda_\mathcal{Y}}{2}\left\Vert \mathcal{Y}_{k+1}-\mathcal{P}\right\Vert_F^2+\lambda_s\left\Vert \mathcal{P} \right\Vert_1.\label{eq P}\\
	\mathcal{Q}_{k+1}&=&\arg\min_{\mathcal{Q}}  \dfrac{\lambda_\mathcal{X}}{2}\left\Vert \mathcal{X}_{k+1}-\mathcal{Q}\right\Vert_F^2.\label{eq Q}
\end{eqnarray}

\begin{itemize}
	\item \textbf{Solving \eqref{eq X} and \eqref{eq Y}:} In certain studies, researchers have proposed to address problems akin to \eqref{eq X} and \eqref{eq Y} utilizing quadratic function tools, as discussed in \cite{imbiriba2019low, imbiriba2}. However, this approach is hindered by its sluggishness.\\
	Recently, the rescaled gradient descent method, also known as the multiplicative updates rule, has yielded promising outcomes for NMF and NTF, as evidenced by \cite{l1_NMF,l12_norm,l2-NMF}. To tackle \eqref{eq X} and \eqref{eq Y}, we will adapt this method to our problem. The following theorem presents the results of the $(k+1)$-th iterations of problems \eqref{eq X} and \eqref{eq Y}.
	\medskip
	\begin{theorem}
		The solutions of \eqref{eq X} and \eqref{eq Y} can be obtained by 
		\begin{eqnarray}
			&\mathcal{X}_{k+1}=\mathcal{X}_{k}.* \left(  \mathcal{A}*_2 \mathcal{Y}_{k}^T  \right)./\left( \mathcal{X}_{k}*_1 \mathcal{Y}_{k}*_2 \mathcal{Y}_{k}^T +\lambda_\mathcal{X}\left( \mathcal{X}_{k}-\mathcal{Q}_{k}  \right)  \right), \label{solX}\\
			&\mathcal{Y}_{k+1}=\mathcal{Y}_{k}.* \left(   \mathcal{X}_{k+1}^T*_1 \mathcal{A}  \right)./ \left( \mathcal{X}_{k+1}^T*_1\mathcal{X}_{k+1}*_1\mathcal{Y}_{k}^T +\lambda_\mathcal{Y}\left(\mathcal{Y}_{k}-\mathcal{P}_{k}  \right)  \right).\label{solY}
		\end{eqnarray}
		where  the symbols $.*$ and $./$ represent, respectively, the element-wise product and division.
	\end{theorem}
	\medskip
	\begin{proof}
		Let's assume the following functions.
		\[
		F(\mathcal{X})=\dfrac{1}{2} \left\Vert \mathcal{A}- \mathcal{X}*_1 \mathcal{Y}_k \right\Vert_F^2 + \dfrac{\lambda_\mathcal{X}}{2}\left\Vert \mathcal{X}-\mathcal{Q}_k\right\Vert_F,
		\]
		and 
		\[
		G(\mathcal{Y})=\dfrac{1}{2} \left\Vert \mathcal{A}- \mathcal{X}_{k+1}*_1 \mathcal{Y} \right\Vert_F^2  + \dfrac{\lambda_\mathcal{Y}}{2}\left\Vert \mathcal{Y}-\mathcal{P}_k\right\Vert_F.
		\]
		Through the utilization of gradient descent methods, we ascertain
		\begin{equation*}
			\mathcal{X}_{k+1}=\mathcal{X}_k - \theta_{\mathcal{X}_k}.* \dfrac{dF(\mathcal{X}_k)}{d\mathcal{X}}, 
		\end{equation*}
		and 
		\begin{equation*}
			\mathcal{Y}_{k+1}=\mathcal{Y}_k - \theta_{\mathcal{Y}_k}.* \dfrac{dG(\mathcal{Y}_k)}{d\mathcal{Y}}.
		\end{equation*}
		Furthermore, 
		\[
		\dfrac{dF(\mathcal{X})}{d\mathcal{X}}=-\left( \mathcal{A}- \mathcal{X}*_1 \mathcal{Y}_{k}  \right)*_2 \mathcal{Y}_k^T + \lambda_\mathcal{X}\left( \mathcal{X}-\mathcal{Q}_k  \right),
		\]
		and
		\[
		\dfrac{dG(\mathcal{Y})}{d\mathcal{Y}}=-\mathcal{X}_{k+1}^T*_1\left( \mathcal{A}- \mathcal{X}_{k+1}*_1 \mathcal{Y}  \right) + \lambda_\mathcal{Y}\left( \mathcal{Y}-\mathcal{P}_k  \right) .
		\]
		If we choose 
		\[
		\theta_{\mathcal{X}_k}= \mathcal{X}_k./\left( \mathcal{X}_k*_1 \mathcal{Y}_k*_1 \mathcal{Y}_k^T  +\lambda_\mathcal{X}\left( \mathcal{X}_k-\mathcal{Q}_k  \right)  \right),\]
		and
		\[
		\theta_{\mathcal{Y}_k}= \mathcal{Y}_k./\left( \mathcal{X}_{k+1}^T*_1 \mathcal{X}_{k+1}*_1 \mathcal{Y}_k^T +\lambda_\mathcal{Y}\left( \mathcal{Y}_k-\mathcal{P}_k  \right)   \right).
		\]
		Subsequently, we obtain the updates provided in Equations \eqref{solX} and \eqref{solY}.
	\end{proof}
	\medskip
	\item \textbf{Solving \eqref{eq P}:} As stated in \cite{prox-l1}, it is evident that the solution to \eqref{eq P} corresponds to the $\ell_1$-norm proximal operator, with the parameter $\dfrac{\lambda_s}{\lambda_\mathcal{Y}}$, since \eqref{eq P} can be expressed as
	\[
	\mathcal{P}_{k+1}=\arg\min_{\mathcal{P}}  \dfrac{1}{2}\left\Vert \mathcal{P}-\mathcal{Y}_{k+1}\right\Vert_F^2+\dfrac{\lambda_s}{\lambda_\mathcal{Y}}\left\Vert \mathcal{P} \right\Vert_1.
	\]
	Thus 
	\begin{equation}\label{solP}
		\mathcal{P}_{k+1}=\left[ \left\vert \mathcal{Y}_{k+1}\right\vert - \dfrac{\lambda_s}{\lambda_\mathcal{Y}} \right].* \texttt{sgn}\left( \mathcal{Y}_{k+1} \right),
	\end{equation}
	where $\left\vert \mathcal{Y}_{k+1}\right\vert$ denotes the tensor of the same size as $\mathcal{Y}_{k+1}$, each of its elements is the absolute value of elements of the tensor $\mathcal{Y}_{k+1}$, and $\left(\texttt{sgn}\left( \mathcal{Y}_{k+1} \right)\right)_{ijk}$ is the sign the $(i,j,k)$-th element of $\mathcal{Y}_{k+1}$.
	\medskip
	\item \textbf{Solving \eqref{eq Q}:} We recall first \eqref{eq Q}; 
	\[
	\mathcal{Q}_{k+1}=\arg\min_{\mathcal{Q}}  \dfrac{\lambda_\mathcal{X}}{2}\left\Vert \mathcal{X}_{k+1}-\mathcal{Q}\right\Vert_F^2.
	\]
	Let us assume that $\mathcal{X}_{k+1}$ exhibits the low-rank property with rank ${\tt rank}_\mathcal{X}$. Consequently, $\mathcal{Q}_{k+1}$ must also be expressed as
	\[
	\mathcal{Q}=\sum_{i=1}^{{\tt rank}_\mathcal{X}} \hat{x}_i^{(1)}\circ \hat{x}_j^{(2)}.
	\]
	Therefore, our problem will be transformed into
	\[
	\left[ \hat{X}^{(1)}, \hat{X}^{(2)} \right]=\arg\min_{ \left[ \hat{X}^{(1)}, \hat{X}^{(2)} \right] }  \dfrac{\lambda_\mathcal{X}}{2}\left\Vert \mathcal{X}_{k+1}-\sum_{i=1}^{{\tt rank}_\mathcal{X}} \hat{x}_i^{(1)}\circ \hat{x}_j^{(2)}\right\Vert_F^2,
	\]
	with $\hat{X}^{(1)}=\left[\hat{x}_1^{(1)},\ldots, \hat{x}_{\texttt{rank}_{\mathcal{X}}}^{(1)}\right]$, and $\hat{X}^{(2)}=\left[\hat{x}_1^{(2)},\ldots, \hat{x}_{\texttt{rank}_{\mathcal{X}}}^{(2)}\right]$. \\
	\begin{equation}\label{solQ}
		\left[ \hat{X}^{(1)}, \hat{X}^{(2)} \right]=\texttt{svds}\left( \mathcal{X}_{k+1}, \texttt{rank}_{\mathcal{X}} \right),
	\end{equation}
	where $\texttt{svds}$ is a Matlab function that returns a subset of singular values and vectors, and $\texttt{rank}_{\mathcal{X}}$ is an estimated rank of $\mathcal{X}$.
\end{itemize}

\medskip
The following algorithm summarizes all the steps of NTF based on the Einstein product.
\begin{algorithm}[H]
	\caption{NTF based on the Einstein product (ENTF).}
	\label{ENTF}
	\textbf{Input:} $\mathcal{A}$, Itermax, $\lambda_s$, $\lambda_\mathcal{X}$, $\lambda_\mathcal{Y}$, $\texttt{rank}_{\mathcal{X}}$.\\
	\textbf{Initialize:} $\mathcal{X}_0$, $\mathcal{Y}_0$, $\mathcal{P}_0$, $\mathcal{Q}_0$.\\
	\textbf{Output:} $\mathcal{X}$ and $\mathcal{Y}$ such that $\mathcal{A} \approx \mathcal{X}*_1 \mathcal{Y}$.\\
	\begin{algorithmic}[1]
		\FOR{$k=1:\text{Itermax}$}
		\STATE Update $\mathcal{X}_{k+1}$ from \eqref{solX}, and take $\mathcal{X}_{k+1}=\max\left\{ 0,\mathcal{X}_{k+1} \right\}$.
		\STATE Update $\mathcal{Y}_{k+1}$ from \eqref{solY}, and take $\mathcal{Y}_{k+1}=\max\left\{ 0,\mathcal{Y}_{k+1} \right\}$.
		\STATE Update $\mathcal{P}_{k+1}$ from \eqref{solP}.
		\STATE Update $\mathcal{Q}_{k+1}$ from \eqref{solQ}.
		\ENDFOR
	\end{algorithmic}
\end{algorithm}

Regarding the stopping criterion employed to assess the convergence of this algorithm, we establish
\[
\dfrac{\left\Vert \mathcal{X}_k - \mathcal{X}_{k-1}\right\Vert_F}{\left\Vert \mathcal{X}_{k-1}\right\Vert_F}\leq \varepsilon \quad \text{and} \quad \dfrac{\left\Vert \mathcal{Y}_k - \mathcal{Y}_{k-1}\right\Vert_F}{\left\Vert \mathcal{Y}_{k-1}\right\Vert_F}\leq \varepsilon,
\]
where $\varepsilon>0$ is a very small positive number.\\

As previously mentioned, in the context of Linear Unmixing, certain constraints are imposed. The first constraint is the abundance non-negativity constraint (ANC), which is incorporated into the NTF Algorithm, alongside the abundance sum-to-one constraint (ASC). To satisfy the latter condition, various methods have been introduced. Some of these methods utilize optimization techniques to ensure compliance \cite{mvnt, imbiriba2019low}, while others employ simple yet effective strategies, as outlined in \cite{augmetation}. In this study, we will adapt the method described in \cite{augmetation} to our problem utilizing the Einstein product. This method involves augmenting the tensors $\mathcal{A}$ and $\mathcal{X}$ so that
\[
\widetilde{\mathcal{A}}=\begin{pmatrix}
	\mathcal{A}\\
	\gamma \textbf{ones}_{J,K}\end{pmatrix}\in \mathbb{R}^{(I+1)\times J\times K},
\]
and 
\[
\widetilde{\mathcal{X}}=\begin{pmatrix}
	\mathcal{X}\\
	\gamma \textbf{1}_r
\end{pmatrix}\in \mathbb{R}^{(I+1)\times r},
\]
where $\gamma>0$ is a parameter to be chosen.

\section{Tensor extrapolation methods}
\label{sec extrapolation}
Extrapolation methods are widely used techniques used to accelerate the convergence of slow algorithms.  These methods consist of transforming the basic sequence into another one that converges faster. First, the extrapolation was used only for scalars \cite{Aitken}, then it was extended to vectors, matrices, and tensors \cite{tahiri,Jbilou1,Jbilou2,Jbilou3,Brezinski,walker2011anderson,Jbilou_extrap1,Jbilou_extrap2,Brezinski2}. They have been applied to many applications, such as the page rank problem \cite{Brezinski} and in image processing, like the completion problem \cite{tahiri}.

Given a sequence of iterates $\mathcal{X}_0,\ldots, \mathcal{X}_n, \ldots $, which are tensors in $\mathbb{R}^{I_1\times \ldots \times I_N}$, assume that this sequence converges to $\widebar{\mathcal{X}}$. The extrapolated sequences are defined as
\[
\mathcal{T}_n^{k}=\sum_{j=0}^{k} \alpha_j \mathcal{X}_{n+j},
\] 
where the coefficients $\alpha_j$ are constrained by the normalization condition $\sum_{j=0}^{k}\alpha_j=1$.\\
This condition is necessary to ensure that when $\left(\mathcal{X}_n\right)_{n\geq 0}$ converges to $\widebar{\mathcal{X}}$, also $\left(\mathcal{T}_n^{k}\right)_{n\geq 0}$ converge to $\widebar{\mathcal{X}}$.\\
To obtain the desired coefficients $\alpha_j$, a projection process such that the following condition is added to the normalization process
\[
\sum_{j=0}^{k}\langle \mathcal{Y}_j, \Delta \mathcal{X}_{n+j} \rangle \alpha_j=0, \; j=0,\ldots, k,
\]
where $\mathcal{Y}_j$ are chosen carefully. The choice of $\mathcal{Y}_j$ determines the method we are going to use, i.e., 
\[
\begin{split}
	\mathcal{Y}_j=\Delta \mathcal{X}_{n+j-1}\; j=1,\ldots, k, \text{ for the MPE method}\\
	\mathcal{Y}_j=\Delta^2 \mathcal{X}_{n+j-1}\; j=1,\ldots, k, \text{ for the RRE method}\\
	\mathcal{Y}_j=\text{arbitrary}\; j=1,\ldots, k, \text{ for the MMPE method}.
\end{split}
\] 
In our case, we are interested in the RRE method with the TET method applied to a series of tensors.\\
In general, both methods consist of producing the following extrapolated sequences
\[
\mathcal{T}_n^k=\mathcal{X}_n - [\Delta \mathcal{X}_n, \Delta \mathcal{X}_{n+1}, \ldots, \Delta \mathcal{X}_{n+k-1} ]\bar{\times}_{N+1} \gamma,
\]
where $[\Delta \mathcal{X}_n, \Delta \mathcal{X}_{n+1}, \ldots, \Delta \mathcal{X}_{n+k-1} ]\in \mathbb{R}^{I_1\times \ldots \times I_N \times k}$ and $\gamma\in \mathbb{R}^k$. The difference between the RRE and the TET methods is about the method of finding the vector $\gamma$.
\begin{itemize}
	\item \textbf{RRE:} The vector $\gamma$ is the solution of the following optimization problem
	\[
	\min_{x} \left\Vert \mathcal{A}\times_{N+1} x - \mathcal{B} \right\Vert_F,
	\] 
	where $\mathcal{A}=[\Delta^2 \mathcal{X}_n,\Delta^2 \mathcal{X}_{n+1}, \ldots, \Delta^2 \mathcal{X}_{n+k-1} ]\in \mathbb{R}^{I_1\times \ldots \times I_N \times k}$ and $\mathcal{B}=\Delta \mathcal{X}_n$.
	\item \textbf{TET:} $\gamma$ is the approximated solution of the problem below
	\[
	\min_{x}  \left\Vert A x - b \right\Vert_F,
	\]
	where $A=\begin{pmatrix}
		\langle \mathcal{Y}, \Delta^2 \mathcal{X}_n \rangle & \ldots & \langle \mathcal{Y}, \Delta^2 \mathcal{X}_{n+k-1} \rangle\\
		\vdots & \ldots & \vdots\\
		\langle \mathcal{Y}, \Delta^2 \mathcal{X}_{n+k-1} \rangle & \ldots & \langle \mathcal{Y}, \Delta^2 \mathcal{X}_{n+2k-2} \rangle 
	\end{pmatrix}$ and $b=\begin{pmatrix}
		\langle \mathcal{Y}, \Delta^2 \mathcal{X}_n \rangle \\
		\vdots\\
		\langle \mathcal{Y}, \Delta^2 \mathcal{X}_{n+k-1} \rangle 
	\end{pmatrix}$, with $\mathcal{Y}\in \mathbb{R}^{I_1\times \ldots \times I_N}$ are carefully chosen.
\end{itemize}
To solve these two least squared problems, we might use the E-QR; Algorithm \ref{alg E-QR}, (or QR) decomposition of $\mathcal{A}$ (or $A$), depending if we are using the RRE or the TET, respectively, i.e.,
\[
\mathcal{A}=\mathcal{Q}\times_{N+1} R^T,\; \text{or}\; A=QR,
\]
with $\mathcal{Q}\in \mathbb{R}^{I_1\times \ldots I_N \times k}$ (or $Q\in \mathbb{R}^{k\times k}$), and $R\in \mathbb{R}^{k\times k}$. In this case
\[
\gamma = R^{-1}y, \; \text{such that}\; y=\mathcal{Q}^T*_{N+1}\mathcal{B},\; \text{or}\; y=Q^Tb.
\]
The following two algorithms sum up the steps of the methods RRE and TTE respectively.
\begin{algorithm}[H]
	\caption{Tensor RRE}
	\label{RRE}
	\textbf{Input:} $\mathcal{X}_n,\ldots, \mathcal{X}_{n+k}$ with $k\geq 2$.\\
	\textbf{Output:} $\mathcal{T}_n^k$.
	\begin{algorithmic}[1]
		\STATE Compute $\Delta \mathcal{X}_i$, $i=n,\ldots, n+k-1$.
		\STATE Compute $\mathcal{A}=[\Delta^2 \mathcal{X}_n, \ldots, \Delta^2 \mathcal{X}_{n+k-2}]$.
		\STATE Compute $\mathcal{A}=\mathcal{Q}*_{N+1} R$ and $\mathcal{B}=\Delta \mathcal{X}_n$.
		\STATE Compute $y=\mathcal{Q}^T*_{N+1}\mathcal{B}$, $\gamma=R^{-1}y$.
		\STATE Compute $\mathcal{T}^{k}_n$.
	\end{algorithmic}
\end{algorithm}
\begin{algorithm}[H]
	\caption{Tensor TET}
	\label{TET}
	\textbf{Input:} $\mathcal{X}_n,\ldots, \mathcal{X}_{n+2k}$ with $k\geq 2$.\\
	\textbf{Output:} $\mathcal{T}_n^k$.
	\begin{algorithmic}[1]
		\STATE Compute $\Delta \mathcal{X}_i$, $i=n,\ldots, n+2k-1$.
		\STATE Compute the matrix $A$ and the vector $b$ as described above.
		\STATE Compute $A=Q R$.
		\STATE Compute $y=Q^T b$, $\gamma=R^{-1}y$.
		\STATE Compute $\mathcal{T}^{k}_n$.
	\end{algorithmic}
\end{algorithm}
These methods of extrapolation will be investigated to be used in the algorithm of the NTF based on the Einstein product (ENTF). The algorithms below describe the ENTF based on the RRE and the TET methods; ENTF-RRE, and ENTF-TET, respectively.
\begin{algorithm}[H]
	\caption{ ENTF-RRE.}
	\label{ENTF-RRE}
	\textbf{Input:} $\mathcal{A}$, Itermax, $\lambda_s$, $\lambda_\mathcal{X}$, $\lambda_\mathcal{Y}$, $\texttt{rank}_{\mathcal{X}}$.\\
	\textbf{Initialize:} $\mathcal{X}_0$, $\mathcal{Y}_0$, $\mathcal{P}_0$, $\mathcal{Q}_0$.\\
	\textbf{Output:} $\mathcal{X}$ and $\mathcal{Y}$ such that $\mathcal{A} \approx \mathcal{X}*_1 \mathcal{Y}$.\\
	\begin{algorithmic}[1]
		\FOR{$k=1:\text{Itermax}$}
		\STATE Compute $\left\{\mathcal{X}_{n}\right\}_{n=1}^k$, $\left\{\mathcal{Y}_{n}\right\}_{n=1}^k$, $\left\{\mathcal{P}_{n}\right\}_{n=1}^k$, $\left\{\mathcal{Q}_{n}\right\}_{n=1}^k$ from Algorithm \ref{ENTF}.
		\STATE Compute $\mathcal{T}_{1, \mathcal{X}}^k$ and $\mathcal{T}_{1, \mathcal{Y}}^k$ associated to the sequences $\left\{\mathcal{X}_{n}\right\}_{n=1}^k$, $\left\{\mathcal{Y}_{n}\right\}_{n=1}^k$ respectively  from Algorithm \ref{RRE}.
		\IF {The algorithm converges}
		\STATE Stop
		\ELSE
		\STATE $\mathcal{X}_1=\mathcal{T}_{1,\mathcal{X}}^k$ and $\mathcal{Y}_1=\mathcal{T}_{1,\mathcal{Y}}^k$.
		\ENDIF
		\ENDFOR
	\end{algorithmic}
\end{algorithm}

\begin{algorithm}[H]
	\caption{ ENTF-TET.}
	\label{ENTF-TET}
	\textbf{Input:} $\mathcal{A}$, Itermax, $\lambda_s$, $\lambda_\mathcal{X}$, $\lambda_\mathcal{Y}$, $\texttt{rank}_{\mathcal{X}}$.\\
	\textbf{Initialize:} $\mathcal{X}_0$, $\mathcal{Y}_0$, $\mathcal{P}_0$, $\mathcal{Q}_0$.\\
	\textbf{Output:} $\mathcal{X}$ and $\mathcal{Y}$ such that $\mathcal{A} \approx \mathcal{X}*_1 \mathcal{Y}$.\\
	\begin{algorithmic}[1]
		\FOR{$k=1:\text{Itermax}$}
		\STATE Compute $\left\{\mathcal{X}_{n}\right\}_{n=1}^{2k}$, $\left\{\mathcal{Y}_{n}\right\}_{n=1}^{2k}$, $\left\{\mathcal{P}_{n}\right\}_{n=1}^{2k}$, $\left\{\mathcal{Q}_{n}\right\}_{n=1}^{2k}$ from Algorithm \ref{ENTF}.
		\STATE Compute $\mathcal{T}_{1, \mathcal{X}}^k$ and $\mathcal{T}_{1, \mathcal{Y}}^k$ associated to the sequences $\left\{\mathcal{X}_{n}\right\}_{n=1}^{2k}$, $\left\{\mathcal{Y}_{n}\right\}_{n=1}^{2k}$ respectively  from Algorithm \ref{TET}.
		\IF {The algorithm converges}
		\STATE Stop
		\ELSE
		\STATE $\mathcal{X}_1=\mathcal{T}_{1,\mathcal{X}}^k$ and $\mathcal{Y}_1=\mathcal{T}_{1,\mathcal{Y}}^k$.
		\ENDIF
		\ENDFOR
	\end{algorithmic}
\end{algorithm}
\section{Numerical experiments}
\label{sec 3}
In this section, experiments are conducted to evaluate Algorithm \ref{ENTF} on Hyperspectral Image (HI) denoising and HI unmixing tasks. The performance of our proposed algorithm will be benchmarked against several established methods commonly used in this domain, including NMF-$\ell_1$, NMF-$\ell_{2}$, NMF-$\ell_{1/2}$, FCLS, SCLSU, MESMA, RUSAL, and MV-NTF-S \cite{l1_NMF,l2-NMF,l12_norm,FCLS,SCLSU,MESMA,RUSAL,mvnt}. The subsequent sections are organized as follows: the first subsection presents the results of HI denoising, while the second subsection is devoted to HI unmixing. All computations are performed on a laptop equipped with 2.3 GHz Intel Core i5 processors and 8 GB of memory using MATLAB 2023b. The initialization of tensors for all experiments is carried out using a fixed seed.\\

The parameters $\lambda_\mathcal{X}$, $\lambda_\mathcal{Y}$, and $\texttt{rank}_{\mathcal{X}}$ are selected heuristically based on the characteristics of the data and its dimensions. While the sparsity parameter is determined to be
\[
\lambda_s=\sum_{i=1}^{I}\dfrac{\left(\sqrt{JK}- \dfrac{\left\Vert \mathcal{A}(i,:,:)\right\Vert_1 }{\left\Vert \mathcal{A}(i,:,:)\right\Vert_F}\right)}{(\sqrt{JK}-1)\sqrt{I}}
\]

To measure the quality of the recovered data, we use:
\begin{itemize}
	\item The Mean Squared Error (MSE) defined by 
	\[
	MSE=\dfrac{\left\Vert \mathcal{A}- \widehat{\mathcal{A}}\right\Vert_F}{IJK},
	\]
	where $\mathcal{A}\in \mathbb{R}^{I\times J\times K}$ is the original HI, and $\widehat{\mathcal{A}}$ is the recovered one.
	\item The Spectral Angle Mapper for the $i$-th endmember
	\[
	SAM_{i}=\arccos\left(\dfrac{ \mathcal{X}(:,i)^T \widehat{\mathcal{X}}(:,i) }{\left\Vert \mathcal{X}(:,i)\right\Vert_F \left\Vert \widehat{\mathcal{X}}(:,i) \right\Vert_F}\right), \quad i=1,\ldots, r,
	\]
	where $\mathcal{X}$ and $\widehat{\mathcal{X}}$ are respectively the original and the obtained endmember, and $r$ is the number of endmembers. In the sequel, $SAM$ denotes the mean of the $SAM$ of the endmembers, i.e., 
	\[
	SAM=\dfrac{1}{r}\sum_{i=1}^{r}SAM_{i}.
	\]
	\item The Mean Squared Error of the abundance
	\[
	MSE_{\mathcal{Y}}=\dfrac{\left\Vert \mathcal{Y}-\widehat{\mathcal{Y}}\right\Vert_F}{IJK},
	\]
	where $\mathcal{Y}$ and $\widehat{\mathcal{Y}}$ denote the original and the obtained abundance tensors.
\end{itemize}
\subsection{Denoising Hyperspectral images}

In this section, we assess the performance of Algorithm \ref{ENTF} in denoising tasks. The datasets employed are sourced from the website \footnote{\url{http://lesun.weebly.com/hyperspectral-data-set.html}}, and they are generated using synthetic tools available in the package \footnote{\url{https://www.ehu.eus/ccwintco/index.php?title=Hyperspectral_Imagery_Synthesis_tools_for_MATLAB}}. These datasets encompass synthetic hyperspectral images (HIs), namely Exponential, Legendre, Matern, Rational, and Spheric. Each HI has spatial dimensions of $128\times 128$ and is composed of $431$ spectral bands. Additionally, four HIs with additive noise have been created from the original synthetic images, resulting in Signal-to-Noise Ratios (SNRs) of $20$, $40$, $60$, and $80$ dB.

\noindent For simplicity, the algorithms are applied exclusively to two HIs: Exponential and Legendre.

\noindent In Table \ref{tab 1}, we depict the Mean Squared Error (MSE), Spectral Angle Mapper (SAM), and MSE$\mathcal{Y}$ acquired from denoising the Exponential and Legendre datasets utilizing NMF-$\ell{1}$, NMF-$\ell_{2}$, NMF-$\ell_{1/2}$, MV-NTF-S, and ENTF for two Signal-to-Noise Ratio (SNR) values: SNR=$20$ and $30$.

\noindent In Figure \ref{fig 1}, the resulting five abundance maps of the Legendre dataset obtained using the methods NMF-$\ell_{1}$, NMF-$\ell_{2}$, NMF-$\ell_{1/2}$, MV-NTF-S, and ENTF are presented. 
\begin{table}[h]
	\centering
	\small\addtolength{\tabcolsep}{-3pt}
	\begin{tabular}{c|c|c|c|c|c|c|c}
		\hline Data &SNR & & NMF-$\ell_1$& NMF-$\ell_2$ & NMF-$\ell_{1/2}$ & MV-NTF-S & ENTF \\
		\hline \multirow{6}{*}{Exponential}& \multirow{3}{*}{20} &MSE&6.97e-05&6.99e-05  &6.33e-05  &1.96e-05 &\textbf{1.79e-05}   \\
		&&SAM&2.20e-01 &2.47e-01  & 9.46e-02 &2.74e-01    &\textbf{5.70e-03}     \\
		&&MSE$_\mathcal{Y}$&1.0e-03& 1.9e-03& 1.1e-03 & 9.19e-04&   \textbf{2.13e-04} \\
		\cline{2-8} &  \multirow{3}{*}{60}&MSE &7.71e-05 &7.73e-05  &6.33e-05  &3.26e-05  &\textbf{3.16e-05}  \\
		& &SAM&2.20e-01& 2.43e-01 &9.46e-02 & 3.304e-01      &   \textbf{1.86e-04}      \\
		&& MSE$_\mathcal{Y}$&1.0e-03 &  1.9e-03&1.1e-03     & 8.91e-04 &\textbf{3.14e-06}\\
		\hline \multirow{6}{*}{Legendre}& \multirow{3}{*}{20} &MSE  & 5.54e-05 & 5.56e-05 &4.80e-05  &1.93e-05  &\textbf{1.92e-05}  \\
		&&SAM& 2.15e-01  & 2.47e-01  &1.21e-01  &4.55e-01     &   \textbf{2.09e-02} \\
		&&MSE$_\mathcal{Y}$& 1.3e-03& 1.9e-03 &1.3e-03  &1.10e-03    & \textbf{1.78e-04}   \\
		\cline{2-8} &  \multirow{3}{*}{60}&MSE  &5.91e-05 &5.93e-05    &  4.80e-05
		&2.59e-05  & \textbf{2.57e-05}\\
		& &SAM&2.13e-01  &2.38e-e1  &1.21e-01   & 4.67e-01      &   \textbf{7.25e-04}    \\
		&& MSE$_\mathcal{Y}$& 1.3e-03& 1.9e-03&1.3e-03    & 1.10e-03 &\textbf{2.53e-06}\\
		\hline
	\end{tabular}
	\caption{The MSE, SAM, and MSE$_\mathcal{Y}$ values obtained through the utilization of the NMF-$\ell_{1}$, NMF-$\ell_{2}$, NMF-$\ell_{1/2}$, MV-NTF-S, and ENTF methods on Exponential and Legendre data with SNR=$20$ and SNR=$60$.}\label{tab 1}
\end{table}

\begin{figure}[h]
	\centering
	\small\addtolength{\tabcolsep}{-5pt}
	\begin{tabular}{c}
		\includegraphics[width=.7\linewidth]{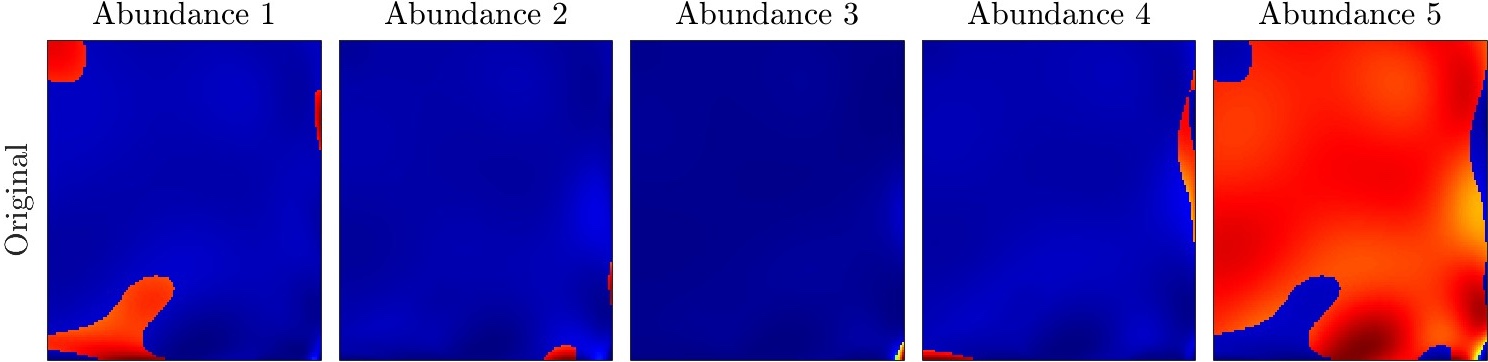} \\
		\includegraphics[width=.7\linewidth]{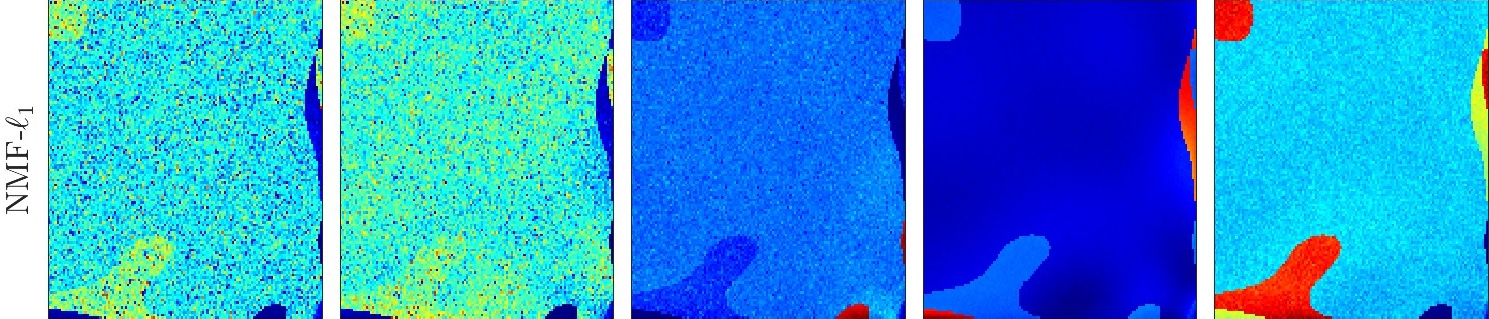} \\
		\includegraphics[width=.7\linewidth]{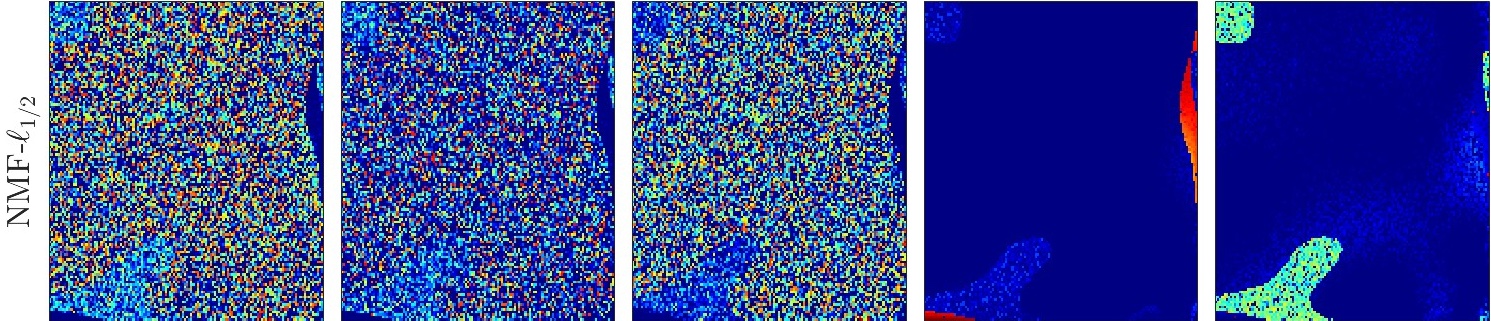} \\
		\includegraphics[width=.7\linewidth]{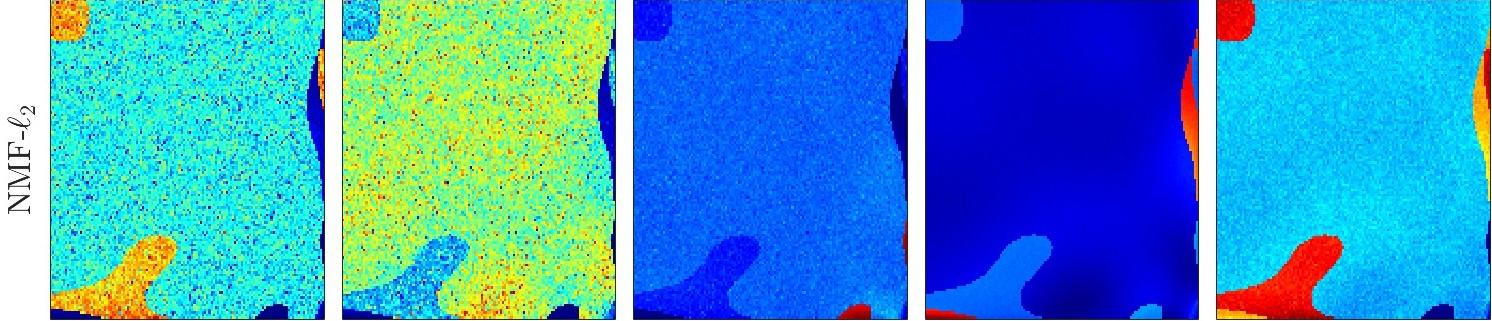}\\
		\includegraphics[width=.7\linewidth]{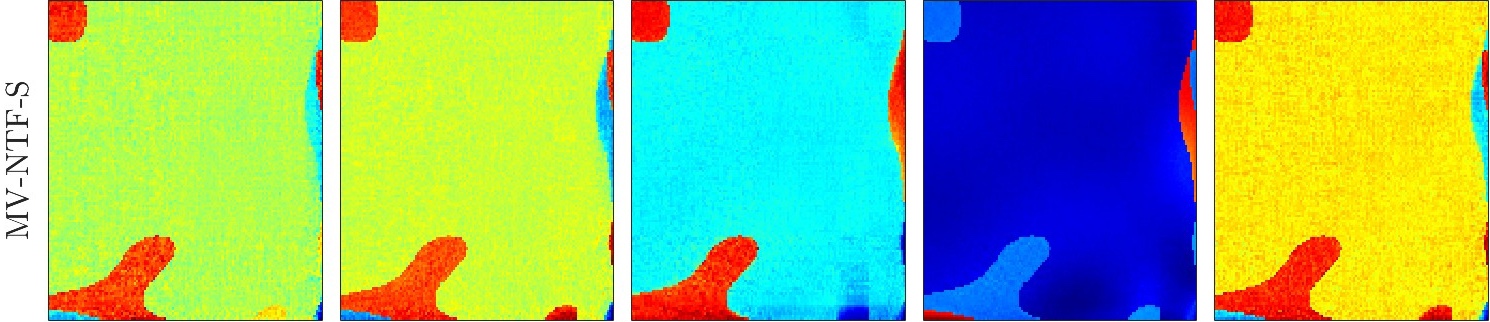}\\
		\includegraphics[width=.7\linewidth]{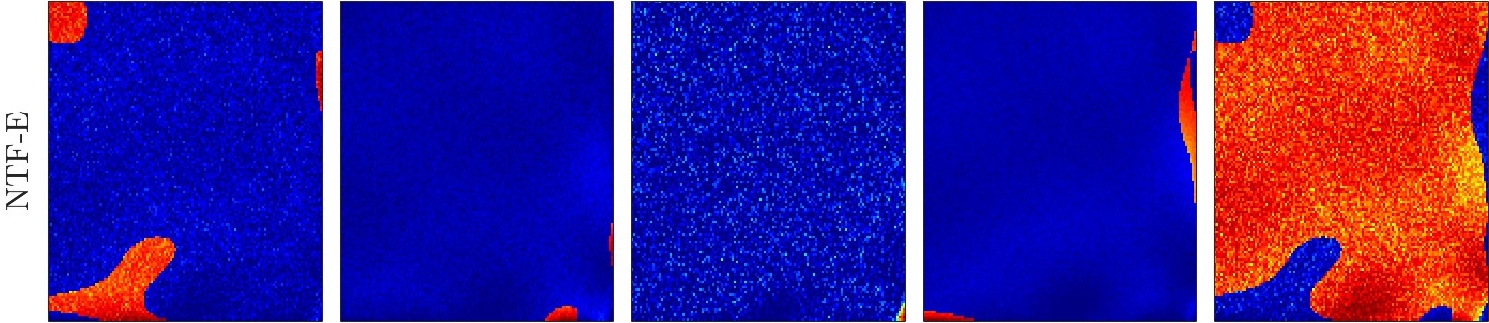}
	\end{tabular}
	\caption{The mapping of the five abundance maps of the Legendre data with SNR=$20$ using the methods NMF-$\ell_{1}$, NMF-$\ell_{2}$, NMF-$\ell_{1/2}$, MV-NTF-S, and ENTF.}\label{fig 1}
\end{figure}

Table \ref{tab 1} and Figure \ref{fig 1} collectively illustrate the effectiveness of the proposed method compared to NMF-$\ell_{1}$, NMF-$\ell_{2}$, NMF-$\ell_{1/2}$, and MV-NTF-S across various Signal-to-Noise Ratio (SNR) levels.
\subsection{Hyperspectral unmixing}

This section is dedicated to Hyperspectral Image Unmixing and is subdivided into two parts: the first part utilizes synthetic data, while the subsequent part employs real data. 
\subsubsection{Synthetic data}
In this section, Hyperspectral Unmixing (HU) is performed on well-established synthetic datasets: DC1, DC2, and DC3, with dimensions $75\times 75\times 224$, $100\times 100\times 224$, and $105\times 105 \times 224$, respectively. The first two indices of these images represent the spatial dimensions, while the last index denotes the spectral bands. It is worth noting that the number of endmembers for DC1, DC2, and DC3 are $5$, $9$, and $6$, respectively.

\begin{figure}[h]
	\centering
	\small\addtolength{\tabcolsep}{-5pt}
	\begin{tabular}{c}
		\includegraphics[width=.7\linewidth]{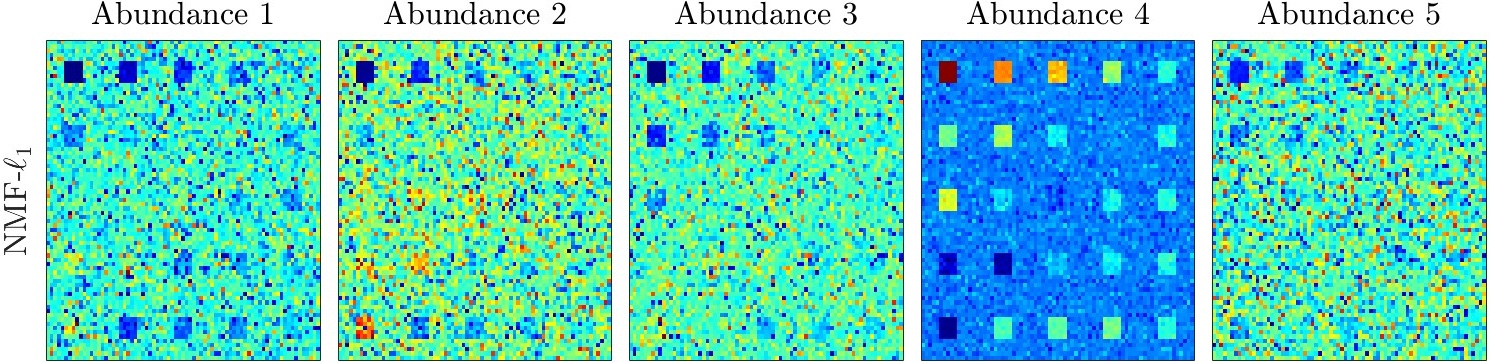} \\
		\includegraphics[width=.7\linewidth]{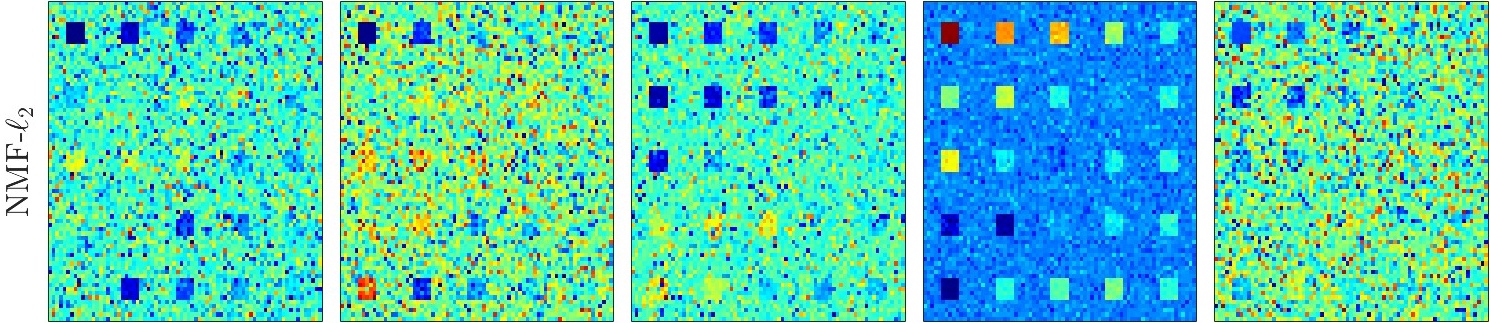} \\
		\includegraphics[width=.7\linewidth]{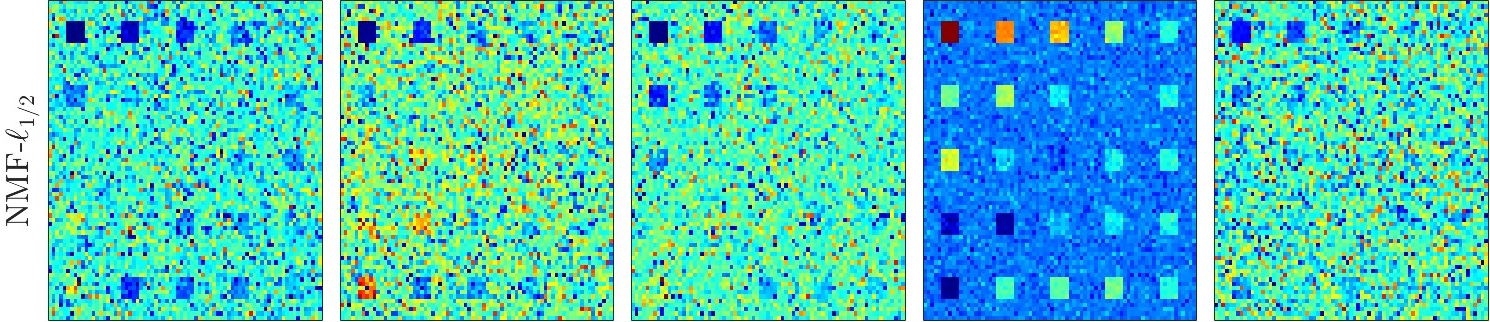}\\
		\includegraphics[width=.7\linewidth]{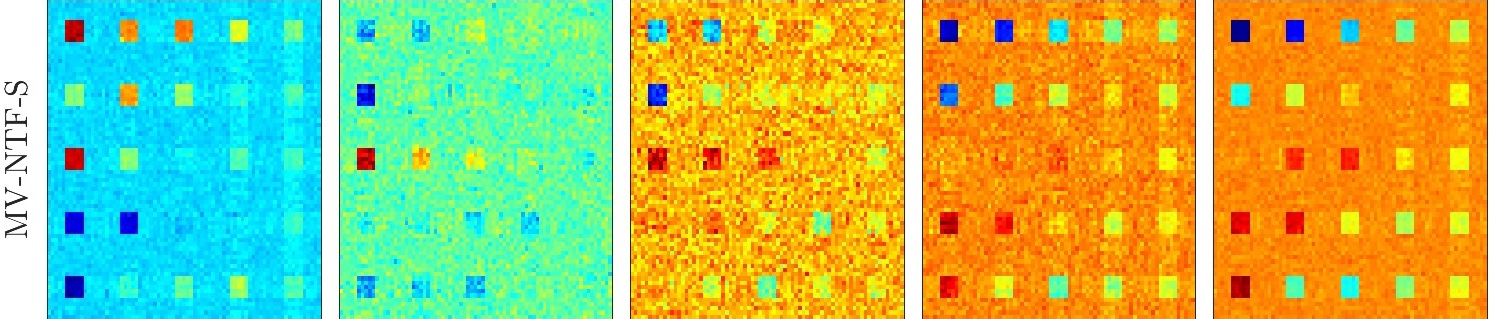}\\
		\includegraphics[width=.7\linewidth]{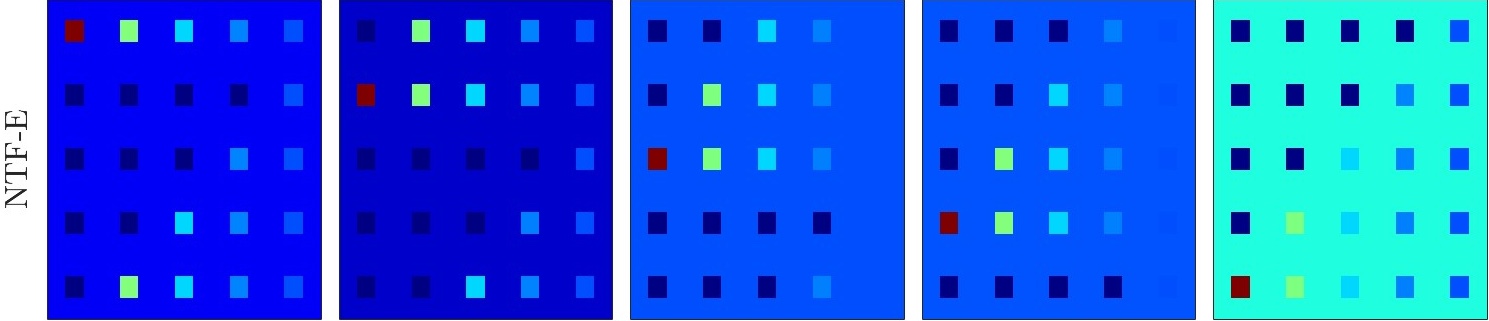}\\
	\end{tabular}
	\caption{The five abundance maps of the DC1 dataset obtained through the utilization of the methods NMF-$\ell_{1}$, NMF-$\ell_{2}$, NMF-$\ell_{1/2}$, MV-NTF-S, and ENTF.}
	\label{fig 2}
\end{figure}

\noindent Figure \ref{fig 2} illustrates the mapping of the five abundances of the DC0 data using the methods NMF-$\ell_{1}$, NMF-$\ell_{2}$, NMF-$\ell_{1/2}$, MV-NTF-S, and ENTF. Additionally, Table \ref{tab 2} presents the values of SAM$i$ with $i=1,\ldots,5$, $i=1,\ldots,9$, $i=1,\ldots,6$, respectively, for the DC1, DC2, and DC3 images, as well as the mean value, obtained using the methods NMF-$\ell{1}$, NMF-$\ell_{2}$, NMF-$\ell_{1/2}$, MV-NTF-S, and ENTF. Furthermore, Table \ref{tab 3} is dedicated to comparing the values of MSE and MSE$_\mathcal{Y}$ obtained through all the aforementioned methods for the DC1, DC2, and DC3 images.

\begin{table}[h]
	\centering
	\small\addtolength{\tabcolsep}{-5pt}
	\begin{tabular}{c|c|c|c|c|c|c}
		\hline Data &  & NMF-$\ell_1$& NMF-$\ell_2$ & NMF-$\ell_{1/2}$ & MV-NTF-S & ENTF \\
		\hline \multirow{6}{*}{DC1}& EM1 & 2.89e-01    &2.88e-01 &2.89e-01 &5.08e-01  & \textbf{5.17e-05}   \\
		& EM2 & 2.00e-01    &1.98e-01  &1.99e-01 &4.82e-01    & \textbf{1.27e-05}  \\
		& EM3 &  2.94e-01   &2.88e-01 & 2.92e-01 &3.94e-01   &  \textbf{5.57e-06}      \\
		& EM4 &2.87e-01    & 2.88e-01  & 2.85e-01 & 4.92e-01  &\textbf{1.17e-05}    \\
		& EM5 &1.09e-01   &1.07e-01 & 1.07e-01   &4.78e-01   &\textbf{9.26e-06} \\
		\cline{2-7} & Mean & 2.36e-01 &2.34e-01 &2.35e-01  & 4.71e-01 & \textbf{1.82e-05} \\
		\hline \multirow{6}{*}{DC2} & EM1 & 4.54e-01 & 4.70e-01 & 1.23e-01 & 4.21e-01 & \textbf{2.22e-06}      \\
		& EM2 &5.11e-01    &3.90e-01  &1.78e-01  &3.72e-01 & \textbf{1.13e-06}       \\
		& EM3 &4.61e-01  &4.81e-01  & 3.37e-01&2.91e-01&   \textbf{1.87e-06}         \\
		& EM4 &1.92e-01  &3.64e-01  &1.91e-01&6.19e-01&  \textbf{1.58e-06}      \\
		& EM5 &3.46e-01 &  3.43e-01 &2.05e-01 & 5.72e-01 &\textbf{1.63e-06}     \\
		& EM6 &3.88e-01 &3.13e-01 & 6.30e-02& 5.41e-01&  \textbf{1.46e-05}     \\
		& EM7 &2.65e-01  & 3.96e-01 & 1.14e-01& 3.77e-01&  \textbf{1.30e-05}      \\
		& EM8 & 2.45e-01 & 2.94e-01 &3.51e-02 & 5.40e-01 &   \textbf{1.60e-06}     \\
		& EM9 & 3.49e-01 & 3.54e-01&1.68e-01& 3.14e-01&    \textbf{1.14e-06}   \\
		\cline{2-7}& Mean &3.57e-01 &3.78e-01& 1.57e-01& 4.49e-01& \textbf{4.31e-06} \\
		\hline \multirow{6}{*}{DC3} &EM1 &1.93e-01  &1.90e-01 &1.89e-01  &5.20e-01   &    \textbf{3.28e-04}          \\
		& EM2 & 1.05e-01   & 8.0e-02&1.29e-01  & 4.73e-01  &  \textbf{1.14e-05}         \\
		& EM3 &2.71e-01  & 2.77e-01 &2.61e-01  & 3.38e-01  &  \textbf{8.36e-05}         \\
		& EM4 & 1.80e-01 & 1.61e-01  &1.90e-01 & 4.18e-01  &  \textbf{5.65e-05}      \\
		& EM5 & 2.82e-01  &2.67e-01 &2.79e-01   &5.41e-01   &    \textbf{3.86e-05}      \\
		& EM6 & 3.06e-01 & 3.15e-01  & 2.68e-01 & 4.68e-01 &    \textbf{8.73e-05}     \\
		\cline{2-7} &Mean & 2.23e-01& 2.15e-01& 2.19e-01& 4.60e-01&  \textbf{1.00e-04}\\
		\hline 
	\end{tabular}
	\caption{The SAM values for each endmember in the DC1, DC2, and DC3 images obtained through the application of the NMF-$\ell_{1}$, NMF-$\ell_{2}$, NMF-$\ell_{1/2}$, MV-NTF-S, and ENTF methods.}\label{tab 2}
\end{table}

\begin{table}[h]
	\centering
	\small\addtolength{\tabcolsep}{-5pt}
	\begin{tabular}{c|c|c|c|c|c|c}
		\hline  Data & \multicolumn{2}{c}{DC1}& \multicolumn{2}{c}{DC2}& \multicolumn{2}{c}{DC3}\\
		\hline Method &  MSE & MSE$_\mathcal{Y}$& MSE & MSE$_\mathcal{Y}$& MSE & MSE$_\mathcal{Y}$\\
		\hline NMF-$\ell_{1}$ & 2.4e-03&8.92e-06 &6.37e-06& 9.03e-04  &  5.21e-06 &5.03e-04  \\
		\hline NMF-$\ell_{2}$  &8.38e-06 & 8.93e-06 &3.87e-06 &1.9e-03 &4.80e-06&1.5e-03 \\
		\hline NMF-$\ell_{1/2}$  &2.1e-03 & 2.0e-03&1.11e-05 &9.27e-04 &5.55e-06 &6.46e-04 \\
		\hline MV-NTF-S  & 1.11e-05&7.42e-04 & 1.30e-05 &9.22e-04 & 9.80e-06& 3.43e-04 \\
		\hline ENTF  & \textbf{2.56e-09}&\textbf{1.07e-07}  &  \textbf{9.66e-10} &  \textbf{6.67e-08} &\textbf{5.31e-09}  & \textbf{1.45e-07}\\
		\hline
	\end{tabular}
	\caption{The MSE and MSE$_\mathcal{Y}$ values for each of the DC1, DC2, and DC3 images obtained through the utilization of the NMF-$\ell_{1}$, NMF-$\ell_{2}$, NMF-$\ell_{1/2}$, MV-NTF-S, and ENTF methods.}\label{tab 3}
\end{table}

Figure \ref{fig 2} demonstrates the superior performance of the proposed method compared to the other methods. Additionally, Table \ref{tab 2} indicates that the ENTF method outperforms NMF-$\ell_{1}$, NMF-$\ell_{2}$, NMF-$\ell_{1/2}$, and MV-NTF-S regardless of the number of endmembers. This is further supported by the low values of MSE and MSE$_\mathcal{Y}$ presented in Table \ref{tab 3}.
\subsubsection{Real data}
In this section, we aim to assess the efficacy of our proposed method using real data. The Hyperspectral images under examination are Jasper Ridge and Samson, both captured by AVIRIS, initially comprising 224 spectral bands. Subsequently, the water absorption bands were removed, resulting in 156 bands for the Samson image and 198 bands for the Jasper Ridge image. The Jasper Ridge dataset has dimensions of $100 \times 100 \times 198$, while the Samson dataset is sized $95 \times 95 \times 156$.

\begin{figure}[h]
	\centering
	\begin{tabular}{c}
		\includegraphics[width=3.5in, height=2.5in]{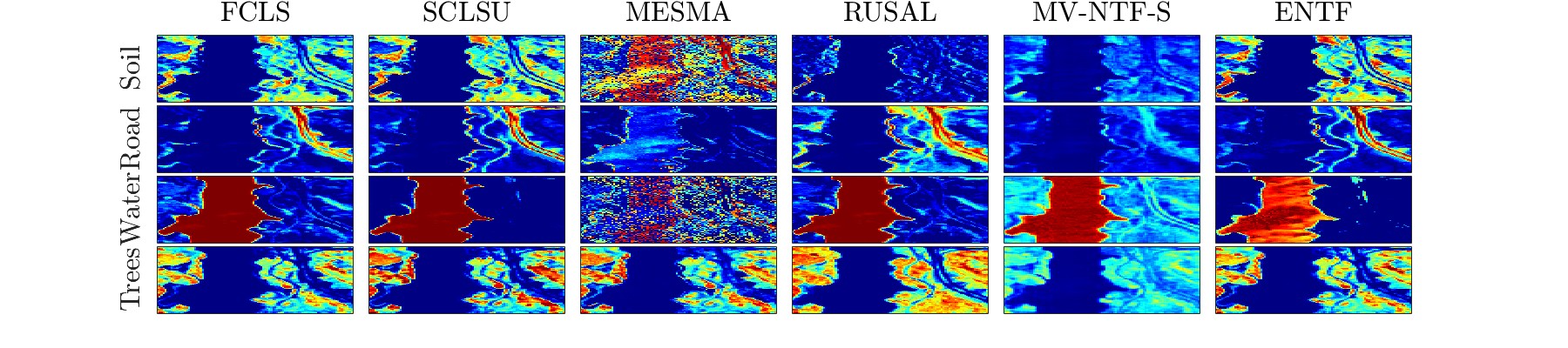}
	\end{tabular}
	\caption{The mapping of the abundances of the endmembers Trees, Soil, Water, and Road for the Jasper Ridge dataset using the FCLS, SCLSU, MESMA, RUSAL, MV-NTF-S, and ENTF methods.}\label{fig 3}
\end{figure}

\begin{figure}[h]
	\centering
	\begin{tabular}{c}
		\includegraphics[width=3.5in, height=2.5in]{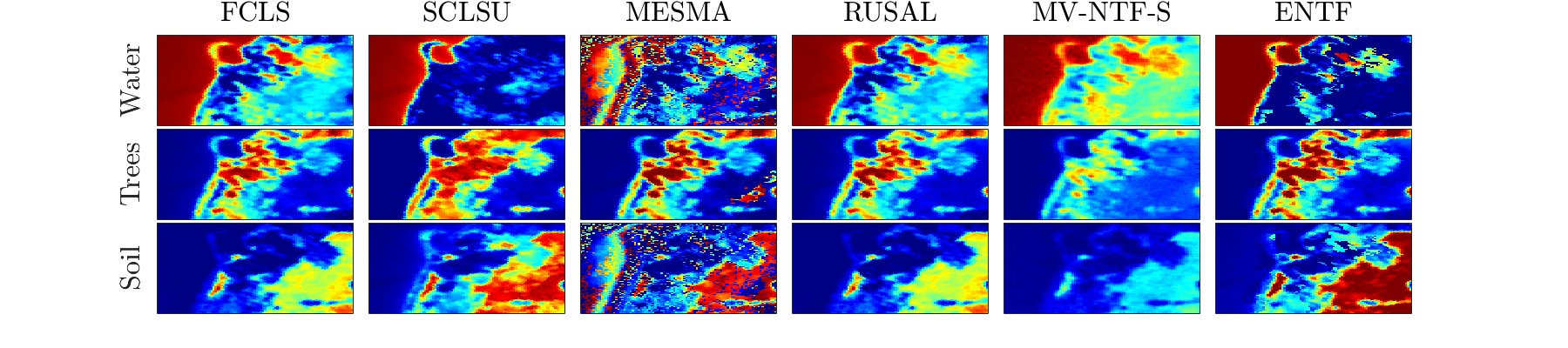}
	\end{tabular}
	\caption{The mapping of the abundances of the endmembers Water, Trees, and Soil for the Samson dataset employing the FCLS, SCLSU, MESMA, RUSAL, MV-NTF-S, and ENTF methods.}\label{fig 4}
\end{figure}

\noindent Table \ref{tab 4} displays the SAM values for each endmember in both the Jasper Ridge and Samson images. Figures \ref{fig 3} and \ref{fig 4} depict the abundance mapping associated with each endmember for the Jasper Ridge and Samson datasets, respectively.
\begin{table}[h]
	\centering
	\begin{tabular}{c|c|c|c|c}
		\hline Data & Endmembers & MESMA & MV-NTF-S & ENTF \\
		\hline\multirow{5}{*}{Jasber Ridge} &Soil & 1.08   &2.24e-01    &    \textbf{1.50e-01}       \\
		& Road& 4.66e-01  & 3.90e-01  &  \textbf{3.34e-02}   \\
		&Water &\textbf{1.70e-01}  & 5.06e-01 & 1.79e-01   \\
		& Trees& 7.95e-02  &  1.265e-01&   \textbf{2.93e-02}   \\
		\cline{2-5} & Mean & 4.50e-01 &3.12e-01 &\textbf{9.82e-02}    \\
		\hline \multirow{4}{*}{Samson}& Soil & 9.51e-01  &2.95e-01  & \textbf{1.18e-01}         \\
		& Trees &7.15e-02 &1.17e-01 & \textbf{7.90e-03}   \\
		& Water &7.52e-02& 1.17e-01&  \textbf{3.11e-02} \\
		\cline{2-5} & Mean &3.66e-01 &1.77e-01   &\textbf{5.24e-02}  \\
		\hline 
	\end{tabular}
	\caption{The SAM values of the endmembers in the Jasper Ridge and Samson datasets were obtained through the utilization of the MESMA, MV-NTF-S, and ENTF methods.}\label{tab 4}
\end{table}

\begin{table}[h]
	\centering
	\small\addtolength{\tabcolsep}{-3pt}
	\begin{tabular}{c|c|c|c|c|c|c}
		\hline Data  &FCLS & SCLSU &MESMA & RUSAL & MV-NTF-S & ENTF \\
		\hline Jasper Ridge& 2.21e-05  & 1.68e-05  & 6.45e-05 & 4.17e-05&  1.23e-05  &  \textbf{9.02e-06}  \\
		\hline Samson & 1.44e-05  &8.44e-06 & 1.57e-04 & 1.60e-05 &8.35e-06 & \textbf{7.25e-06}   \\
		\hline 
	\end{tabular}
	\caption{The MSE values obtained by using the methods MESMA, MV-NTF-S, and ENTF, for Jasper Ridge and Samson datasets.}\label{tab 444}
\end{table}

Tables \ref{tab 4} and \ref{tab 444} demonstrate the effectiveness of the proposed method (ENTF) compared to the other methods employed. It is noteworthy that the method performs less effectively only for the endmember associated with Water in the Jasper Ridge dataset.

Figures \ref{fig 3} and \ref{fig 4} demonstrate that the proposed method effectively identifies the desired materials in both the Jasper Ridge and Samson datasets.
\begin{figure}[h]
	\centering
	\subfloat[Soil]{\includegraphics[width=.49\linewidth]{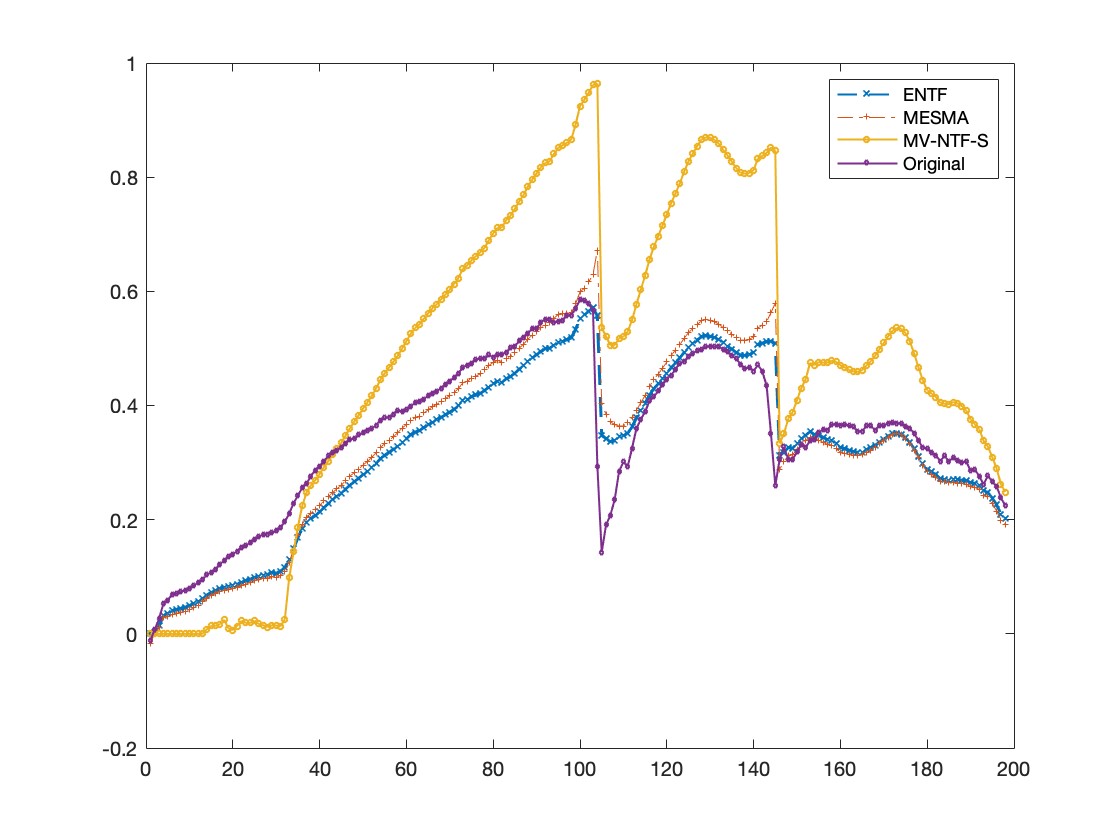}}\hfill
	\subfloat[Road]{\includegraphics[width=.49\linewidth]{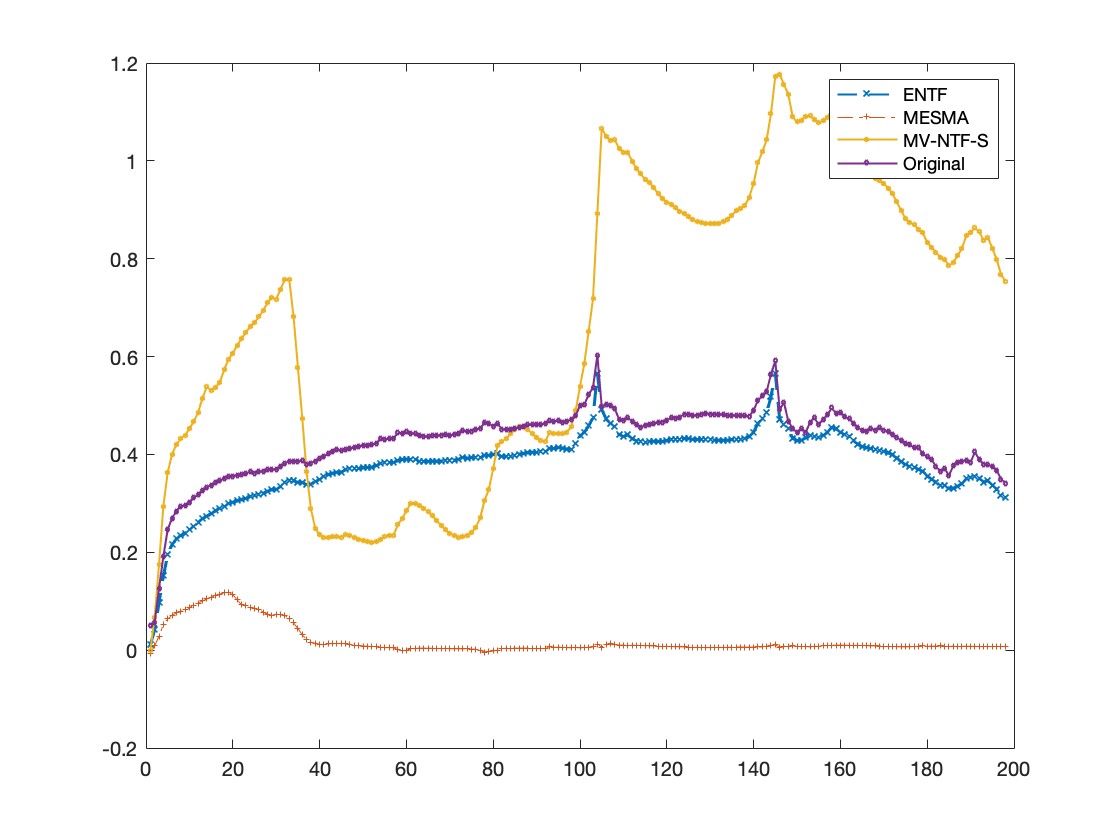}}\hfill
	\subfloat[Water]{\includegraphics[width=.49\linewidth]{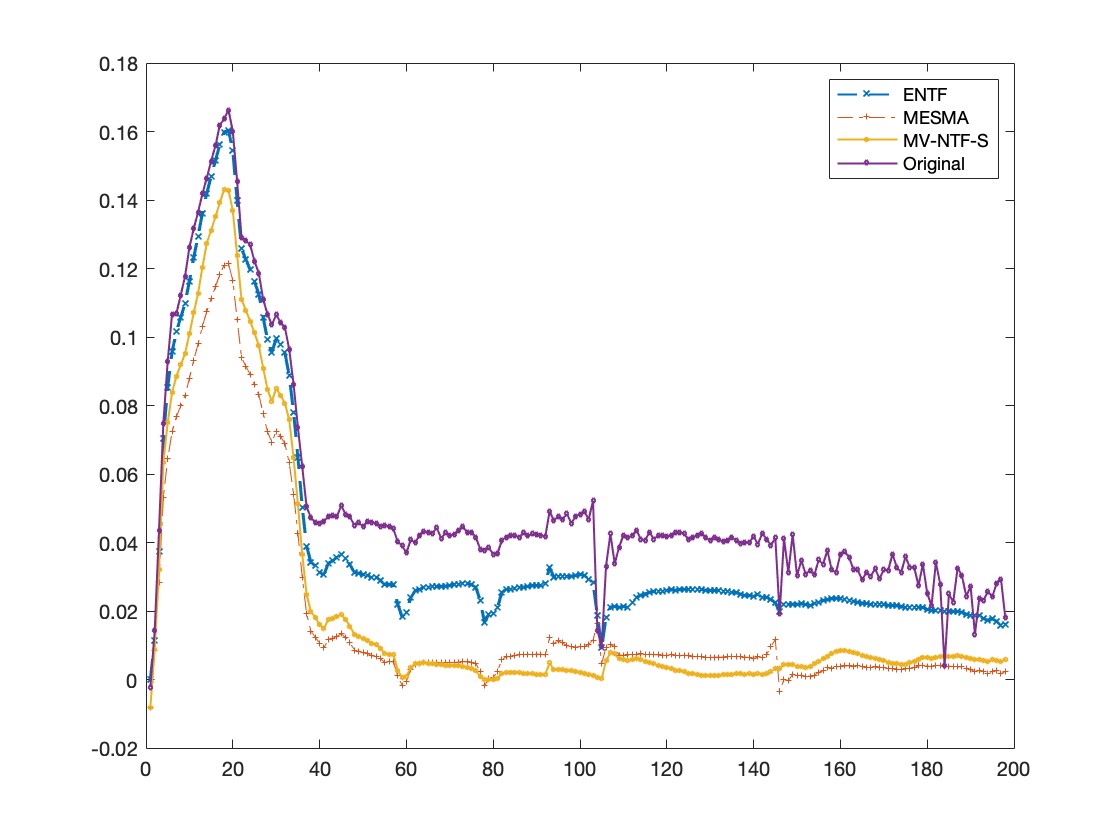}}\hfill
	\subfloat[Trees]{\includegraphics[width=.49\linewidth]{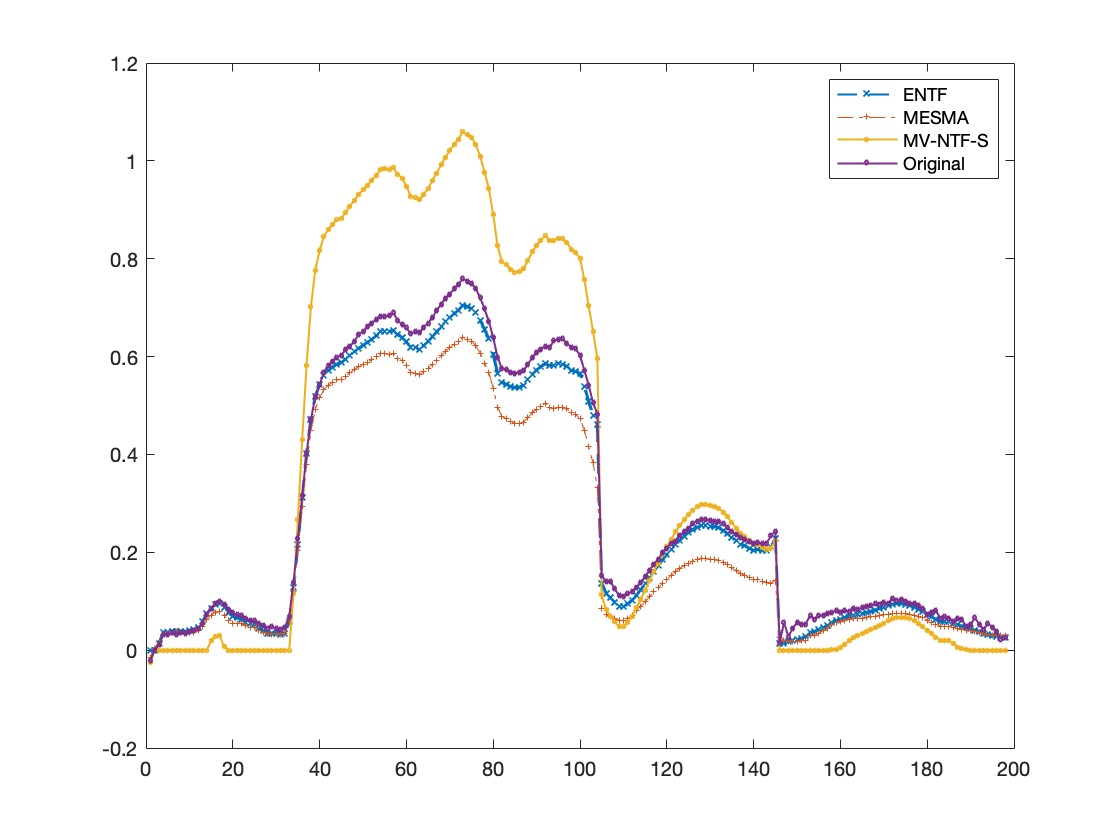}}
	\caption{Jasper Ridge endmembers obtained by the methods MESMA, MV-NTF-S, and ENTF.}\label{fig 33}
\end{figure}
\begin{figure}[h]
	\centering
	\subfloat[Water]{\includegraphics[width=.49\linewidth]{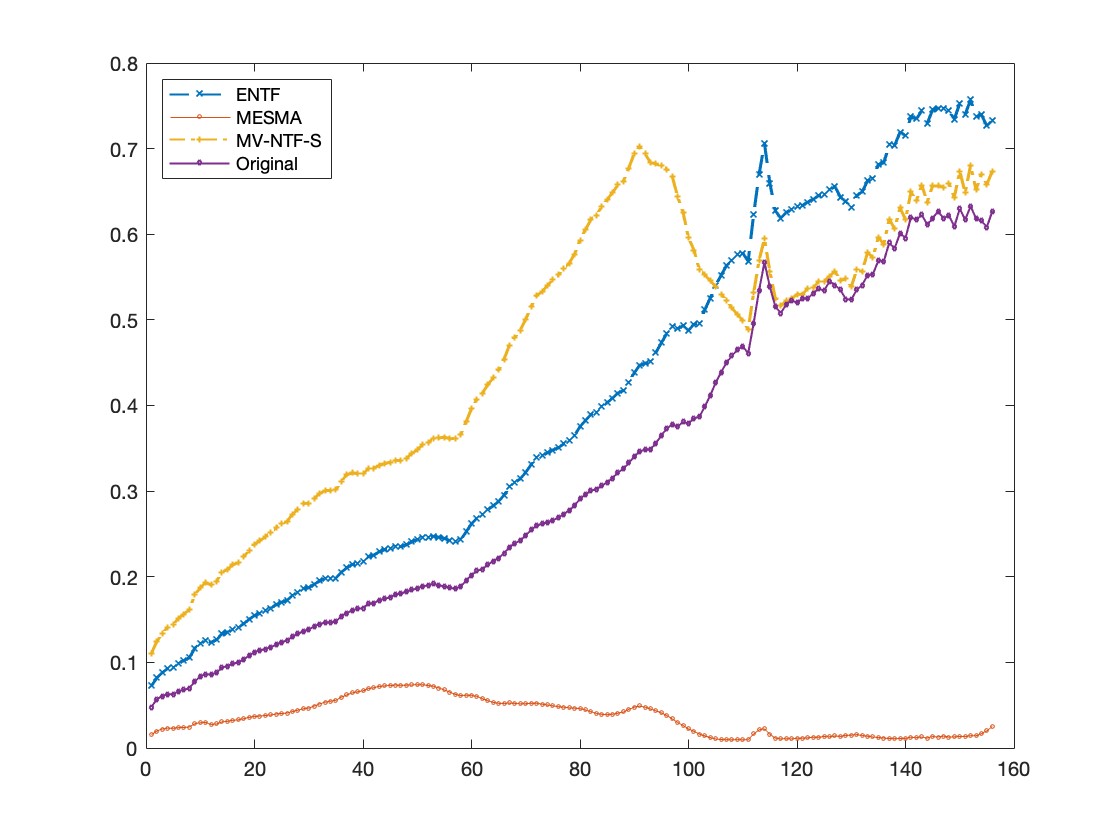}}\hfill
	\subfloat[Trees]{\includegraphics[width=.49\linewidth]{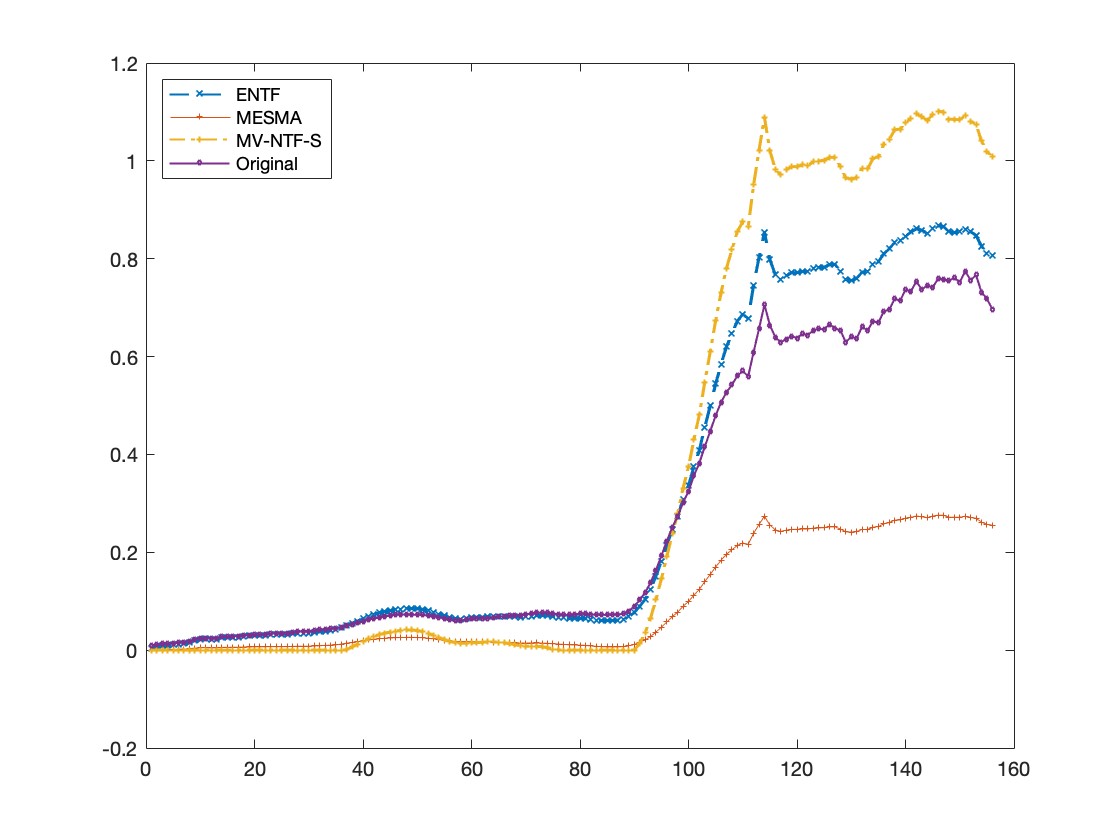}}\hfill
	\subfloat[Soil]{\includegraphics[width=.49\linewidth]{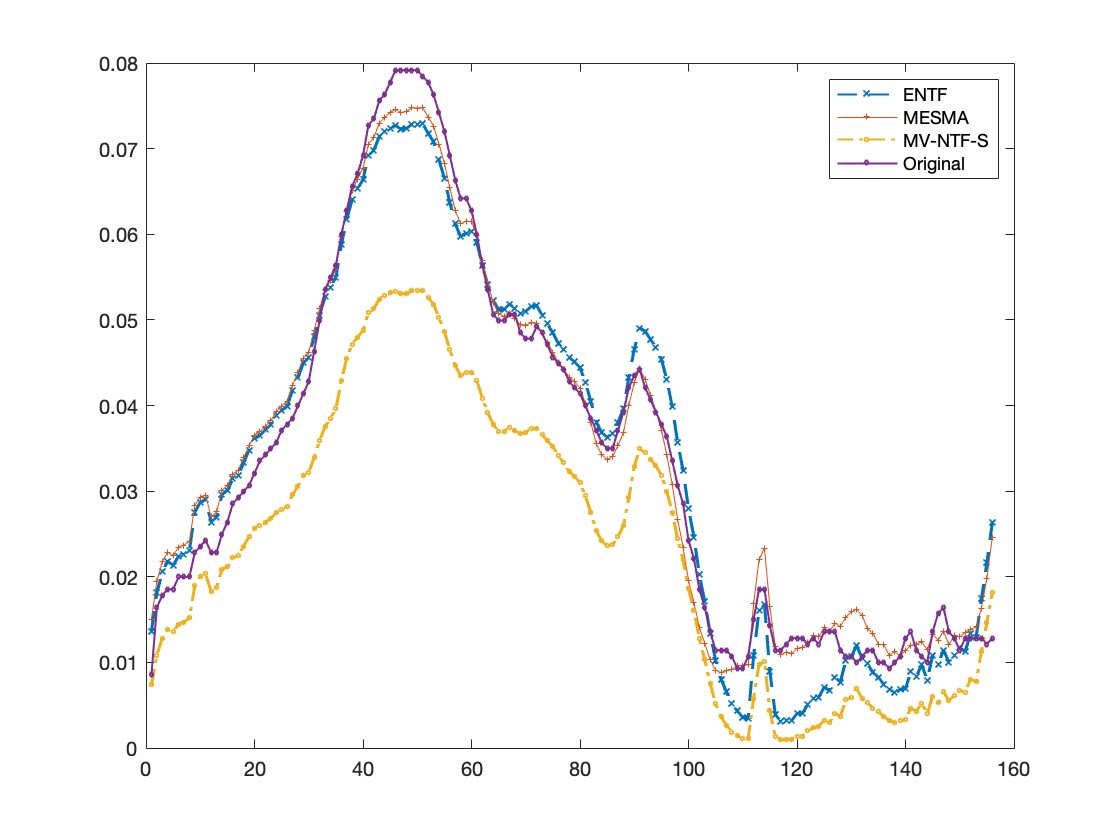}}
	\caption{Samson's endmembers obtained by the methods MESMA, MV-NTF-S, and ENTF.}\label{fig 44}
\end{figure}

Figures \ref{fig 33} and \ref{fig 44} illustrate that the ENTF method produces endmembers closest to the original ones compared to the other methods.
\subsection{Results of the acceleration methods}

In this section, we will compare the proposed ENTF method with methods utilizing extrapolation techniques, namely, ENTF-RRE and ENTF-TET. Specifically, we will focus on the Jasper Ridge dataset for this comparison.

\begin{figure}[h]
	\centering
	\subfloat{\includegraphics[width=.49\linewidth]{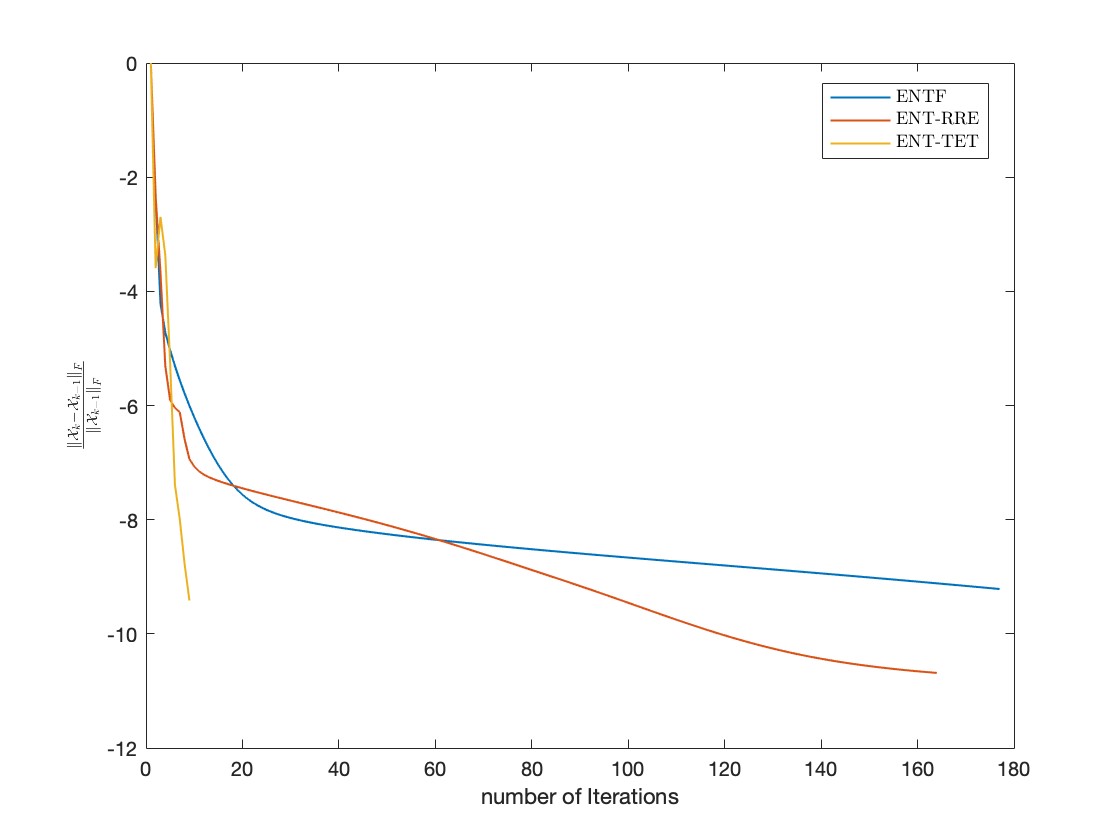}} \hfill
	\subfloat{\includegraphics[width=.49\linewidth]{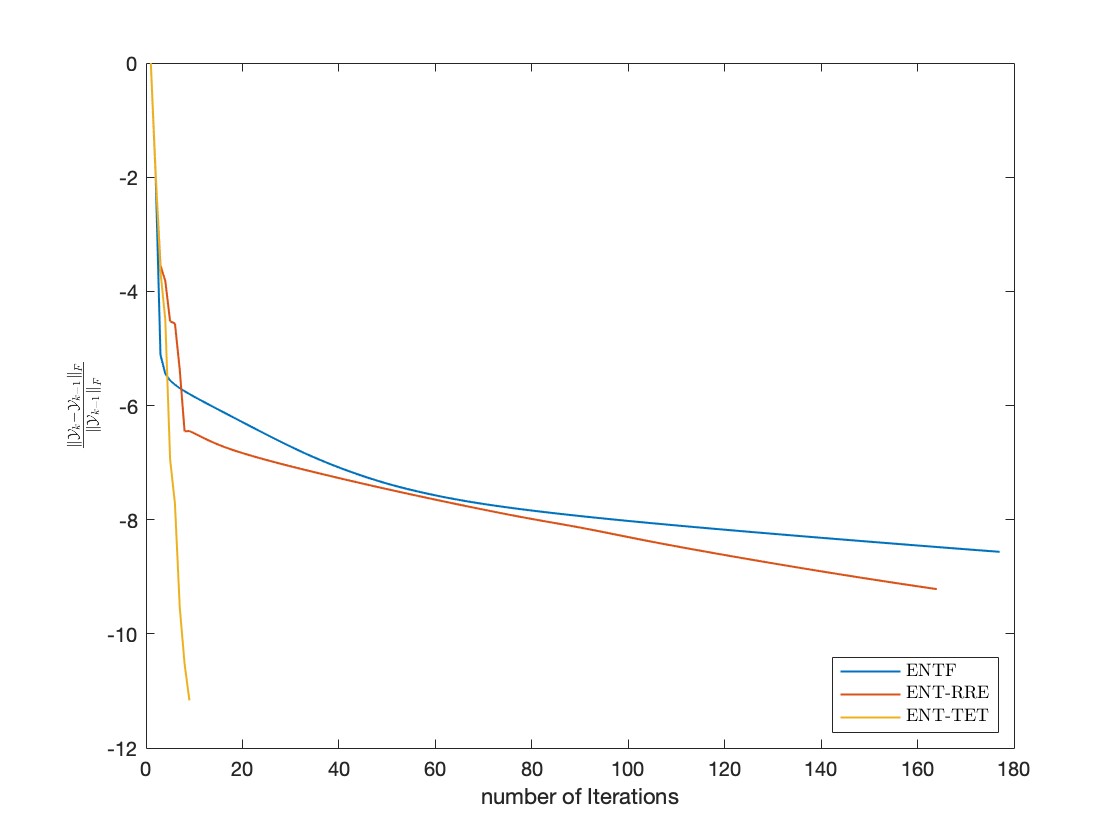}}
	\caption{The errors $\dfrac{\left\Vert \mathcal{X}_k - \mathcal{X}_{k-1}\right\Vert_F}{\left\Vert \mathcal{X}_{k-1}\right\Vert_F}$ and $\dfrac{\left\Vert \mathcal{Y}_k - \mathcal{Y}_{k-1}\right\Vert_F}{\left\Vert \mathcal{Y}_{k-1}\right\Vert_F}$, respectively from the left to the right by using the methods ENTF, ENTF-RRE, and ENTF-TET on the Jasper Ridge dataset.}\label{fig 67}
\end{figure}

\begin{figure}[h]
	\centering
	\subfloat[Soil]{\includegraphics[width=.49\linewidth]{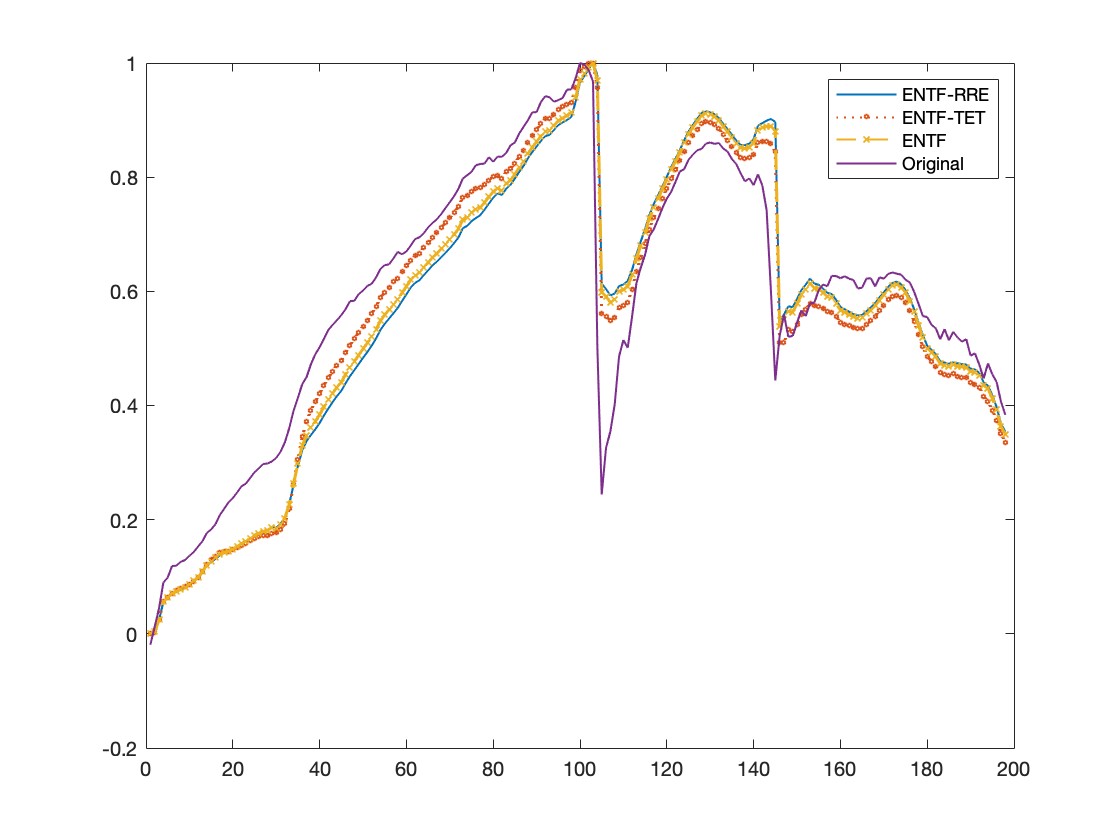}} \hfill
	\subfloat[Road]{\includegraphics[width=.49\linewidth]{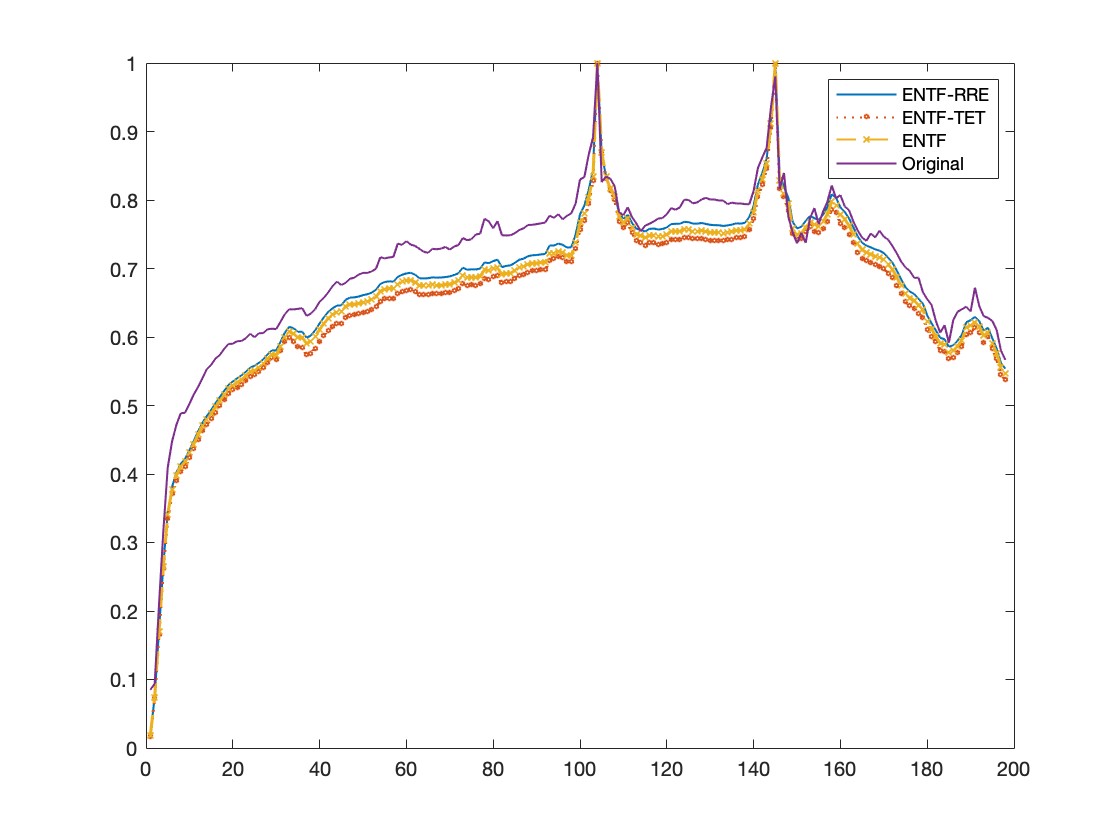}} \hfill
	\subfloat[Water]{\includegraphics[width=.49\linewidth]{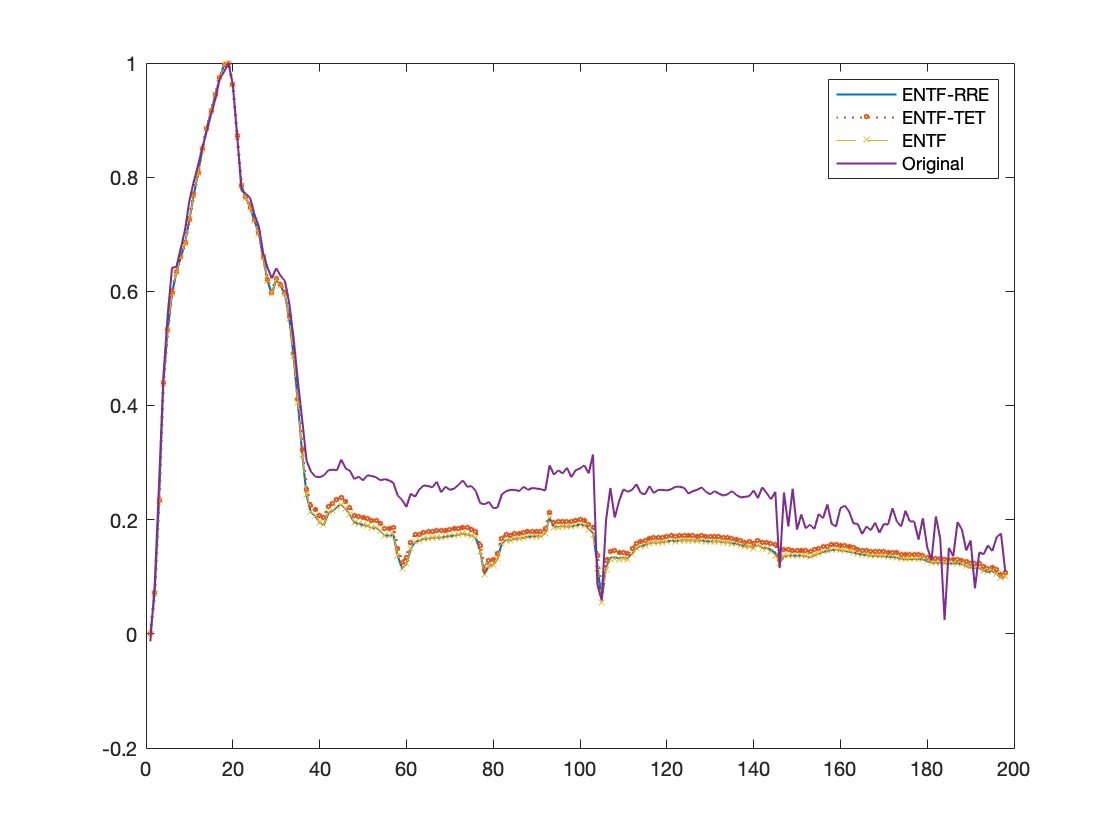}} \hfill
	\subfloat[Trees]{\includegraphics[width=.49\linewidth]{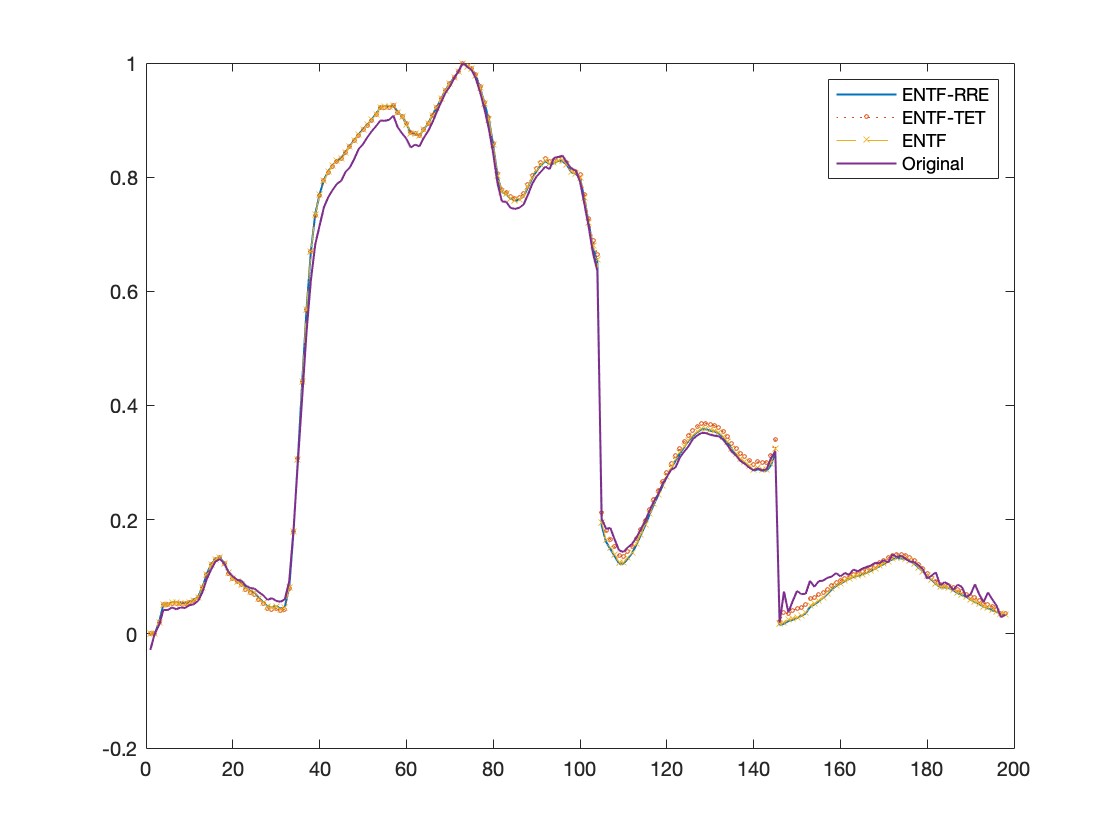}}
	\caption{Jasper Ridge endmembers obtained by the methods ENTF, ENTF-RRE, and ENTF-TE.}\label{fig}
\end{figure}

\noindent Table \ref{tab 66} presents the values of the Mean Squared Error (MSE), Spectral Angle Mapper (SAM), and the number of iterations (iter) obtained by employing the methods ENTF, ENTF-RRE, and ENTF-TET. Additionally, Figure \ref{fig 67} illustrates the error of the sequences $\left\{ \mathcal{X}_k\right\}{k\geq 0}$ and $\left\{ \mathcal{Y}_k\right\}{k\geq 0}$ of the endmembers and the abundance, respectively, at each iteration. Furthermore, the endmembers obtained by the methods ENTF, ENTF-RRE, and ENTF-TET are plotted in Figure \ref{fig}.

\begin{table}[h]
	\centering
	\begin{tabular}{c|c|c|c}
		\hline Method & MSE & SAM & Iter  \\
		\hline ENTF &4.92e-02 & 9.71e-02& 177 \\
		\hline ENTF-RRE &4.89e-02 &9.94e-02 &164   \\
		\hline ENTF-TET & 5.19e-02&8.73e-02 & 9 \\
	\end{tabular}
	\caption{The values of MSE and the SAM and the number of iterations (Iter) of the methods ENTF, ENTF-RRE, and ENTF-TET on the Jasper Ridge dataset.} \label{tab 66}
\end{table}

Based on Table \ref{tab 66} and Figure \ref{fig}, it can be observed that the extrapolation methods yield results comparable to those obtained without extrapolation. Additionally, as depicted in Table \ref{tab 66} and Figure \ref{fig 67}, the sequences $\left\{\mathcal{X}_k\right\}{k\geq 0}$ and $\left\{\mathcal{Y}_k\right\}{k\geq 0}$ converge more rapidly when extrapolation methods are employed, requiring fewer iterations.
\section{Conclusion}
This paper introduces a novel method for Non-Negative Tensor Factorization (NTF) utilizing the Einstein product and incorporating regularization techniques to maintain both smoothness and sparsity within the data. The primary optimization problem is addressed using the adapted Multiplicative Updates method. The methods of extrapolation RRE and TET are adapted to this method. The efficacy of this approach is demonstrated through numerical experiments conducted for denoising and Hyperspectral Image Unmixing, and the extrapolated methods provide a good impression.


\begin{thebibliography}{99}
	\bibitem{cproduct}
	S. Aeron, E. Kernfeld, M. Kilmer, Tensor-tensor products with 
	invertible linear transforms, Linear Algebra and its Applications, 485, 545–570, (2015).
	
	\bibitem{Aitken}
	A. C. Aitken, On Bernoulli’s numerical solution of algebraic equations,
	Proceedings of the Royal Society of Edinburgh 46, 289–305, (1925).
	
	\bibitem{supervised}
	Y. Altmann, N. Dobigeon, A. Halimi,  J. Y. Tourneret, Supervised nonlinear spectral 
	unmixing using a postnonlinear mixing model for hyperspectral imagery, IEEE Transactions on Image Processing, 21(6), 3017--3025, (2012).
	
	\bibitem{unsupervised1}
	Y. Altmann, N. Dobigeon, J. Y. Tourneret, Unsupervised post-nonlinear unmixing of 
	hyperspectral images using a Hamiltonian Monte Carlo algorithm, IEEE Transactions on Image Processing, 23(6), 2663--2675, (2014).
	
	\bibitem{book1}
	K. V.  Arya, N. U. Khan, S. S. Rajput,  A. K. Singh,  Digital Image Enhancement and 
	Reconstruction, (2022).
	
	\bibitem{l0-NP_hard2}
	M. Avellaneda, G. Davis, S. Mallat, Adaptive greedy approximations,
	Constructive Approximation, 13(1), 57–-98, (1997).
	
	\bibitem{kolda}
	B. W. Bader, T. Kolda. Tensor decompositions and applications,
	SIAM review, 51(3), 455–500, (2009).
	
	\bibitem{prox-l1}
	A. Beck, First-order methods in optimization, SIAM, (2017).
	
	\bibitem{Jbilou_extrap1}
	F. P. A. Beik, A. El Ichi, K. Jbilou, R. Sadaka, Tensor extrapolation methods 
	with applications, Numerical Algorithms, 87, 1421–1444, (2021).
	
	\bibitem{elha1}
	A. H. Bentbib, A. El Hachimi, K. Jbilou, A. Ratnani, A tensor 
	regularized nuclear norm method for image and video completion, Journal of Optimization Theory and Applications, 192(2), 401--425, (2022).
	
	\bibitem{elha2}
	A. H. Bentbib, A. El Hachimi, K. Jbilou, A. Ratnani, Fast 
	multidimensional completion and principal component analysis methods via the cosine product, Calcolo, 59(3), 26, (2022).
	
	\bibitem{tahiri}
	A. H. Bentbib, K. Jbilou, R. Tahiri, Multidimensional extrapolated 
	global proximal gradient and applications for image processing, arXiv preprint arXiv:2401.03031, (2024).  
	
	\bibitem{imbiriba2}
	J. C. M Bermudez, R. A. Borsoi,   T. Imbiriba,  A low-rank tensor 
	regularization strategy for hyperspectral unmixing, 2018 IEEE Statistical Signal Processing Workshop (SSP), 373--377, (2018).
	
	\bibitem{imbiriba2019low}
	J. C. M Bermudez,  R. A. Borsoi, T. Imbiriba , Low-rank tensor 
	modeling for hyperspectral unmixing accounting for spectral variability, IEEE Transactions on Geoscience and Remote Sensing, 58(3), 1833-1842, (2019).
	
	\bibitem{HI1}
	J.M. Bioucas-Dias, J. Chanussot, N. Dobigeon, Q. Du,   P. Gader, M. Parente, A. Plaza,   
	Hyperspectral unmixing overview: Geometrical, statistical, and sparse regression-based approaches, IEEE journal of selected topics in applied earth observations and remote sensing, 5(2), 354-379, (2012).
	
	\bibitem{RUSAL}
	J. M. Bioucas-Dias, G.S. Buller, N. Dobigeon, A. Halimi,   S. McLaughlin, Fast 
	hyperspectral unmixing in presence of nonlinearity or mismodeling effects, IEEE Transactions on Computational Imaging, 3(2), 146-159, (2016). 
	
	\bibitem{einstein2}
	M. Brazell, N. Li, C. Navasca, C. Tamon, Solving multilinear 
	systems via tensor inversion, SIAM Journal on Matrix Analysis and Applications, 34(2), 542–570, (2013).
	
	\bibitem{Brezinski2}
	C. Brezinski, Généralisation de la transformation de Shanks, de la table de 
	Padé et l’epsilon algorithm, Calcolo, 12, 317–360, (1975).
	
	\bibitem{Brezinski}
	C. Brezinski, M. Redivo-Zaglia, S. Serra-Capizzano, Extrapolation methods for 
	pagerank computations, Comptes Rendus, Math{\'e}matique, 340(5), 393–397, (2005).
	
	
	\bibitem{einstein3} 
	C. Bu, L. Sun, Y. Wei, B. Zheng, Moore–Penrose 
	inverse of tensors via Einstein product, Linear and Multilinear Algebra, 64(4), 686–698, (2016).
	\bibitem{augmetation}
	C.-I. Chang, D. C. Heinz,  Fully constrained least squares linear spectral 
	mixture analysis method for material quantification in hyperspectral imagery, IEEE Transactions on Geoscience and Remote Sensing, 39(3), 529–545, (2001).
	
	\bibitem{semisupervised}
	C. I. Chang, N. Dobigeon, J. Y. Tourneret,  Semi-supervised linear spectral unmixing 
	using a hierarchical Bayesian model for hyperspectral imagery, IEEE Transactions on Signal Processing, 56(7), 2684--2695, (2008).
	
	\bibitem{SCLSU}
	J. Chanussot, L. Drumetz, S. Henrot, C. Jutten, R. Phlypo, M.-A. Veganzones, Blind hyperspectral unmixing using an extended linear
	mixing model to address spectral variability, IEEE Transactions on 
	Image Processing, 25(8), 3890–3905, (2016).
	
	\bibitem{unsupervised2}
	M. Chen, S. Godsill, J. Y. Tourneret, Q. Wei,  Unsupervised nonlinear spectral unmixing 
	based on a multilinear mixing model, IEEE Transactions on Geoscience and Remote Sensing, 55(8), 4534--4544, (2017).
	
	\bibitem{MESMA}
	R. Church, M. Gardner, R.O. Green, D.A. Roberts,  G. Scheer, S. Ustin,  Mapping 
	chaparral in the Santa Monica Mountains using multiple endmember spectral mixture models, Remote sensing of environment, 65(3), 267-279, (1998).
	
	\bibitem{VCA}
	J. M. Dias, J. M. Nascimento,  Vertex component analysis: A fast algorithm to unmix 
	hyperspectral data, IEEE transactions on Geoscience and Remote Sensing, 43(4), 898--910, (2005).
	
	\bibitem{einstein1}
	A. Einstein, The collected papers of Albert Einstein, 14, (1987).
	
	\bibitem{elha8}
	A. El Hachimi, K. Jbilou, M. Hached, A. Ratnani, Tensor Golub 
	Kahan based on Einstein product, arXiv preprint arXiv:2311.03109, (2023). 
	
	\bibitem{elhalouy}
	S. El-Halouy, S. Noschese, L. Reichel, Perron communicability and 
	sensitivity of multilayer networks, Numerical Algorithms, 92, 597–617, (2023).
	
	\bibitem{Jbilou_extrap2}
	A. El Ichi, K. Jbilou, R. Sadaka, Tensor Global Extrapolation Methods Using 
	the n-Mode and the Einstein Products, Mathematics, 8(8), 1298, (2020).
	
	\bibitem{BSS}
	P. Garat, D. T. Pham,  Blind separation of mixture of independent sources through a 
	quasi-maximum likelihood approach, IEEE transactions on Signal Processing, 45(7), 1712--1725, (1997).
	
	\bibitem{gillis}
	N. Gillis, non-negative matrix factorization, SIAM, (2020).
	
	\bibitem{wu2023multidimensional}
	M. S. Guo, J. Huang,  L. Wu,  Multidimensional 
	low-rank representation for sparse hyperspectral unmixing, IEEE Geoscience and Remote Sensing Letters, 20, 1--5, (2023).
	\bibitem{FCLS}
	D. C. Heinz, Fully constrained least squares linear spectral mixture analysis 
	method for material quantification in hyperspectral imagery, IEEE transactions on geoscience and remote sensing, 39(3), 529-545, (2001).
	
	\bibitem{l1_NMF}
	P. O. Hoyer, Non-negative sparse coding, Proceedings of the 12th IEEE 
	workshop on neural networks for signal processing, 557--565, (2002).
	
	\bibitem{ICA}
	A. Hyv{\"a}rinen, E. Oja, Independent component analysis: algorithms and 
	applications, Neural networks, 13(4-5), 411--430, (2000).
	
	
	\bibitem{Jbilou1}
	K. Jbilou, H. Sadok, Matrix polynomial and epsilon-type extrapolation methods 
	with applications, Numerical Algorithms, 68, 107–119, (2015).
	
	\bibitem{Jbilou2}
	K. Jbilou, H. Sadok, Vector extrapolation methods. applications and numerical 
	comparison, Journal of computational and applied mathematics, 122(1-2), 149–165, (2000).
	
	\bibitem{Jbilou3}
	K. Jbilou, A. Messaoudi, Block extrapolation methods with applications,
	Applied Numerical Mathematics, 106, 154–164, (2016).
	
	\bibitem{l12_norm}
	S. Jia, Y. Qian, A. Robles-Kelly,  J. Zhou,   Hyperspectral 
	unmixing via $ L\_\{1/2\}$ 
	sparsity-constrained non-negative matrix factorization, IEEE Transactions on Geoscience and Remote Sensing, 49(11), 4282--4297, (2011).
	
	\bibitem{tproduct}
	M. E. Kilmer, C. D. Martin, Factorization strategies for third- 
	order tensors, Linear Algebra and its Applications, 435, 641-658, (2011).
	
	\bibitem{book_medecin}
	R. Koprowski, Processing of Hyperspectral Medical Images, (2017).
	
	\bibitem{Multi_updaes}
	D. Lee, H. S, Seung, Algorithms for non-negative matrix factorization,
	Advances in neural information processing systems, 13, (2000).
	
	\bibitem{walker2011anderson}
	P. Ni, H. F. Walker,  Anderson acceleration for fixed-point iterations,
	SIAM Journal on Numerical Analysis, 49(4), 1715--1735, (2011).
	
	\bibitem{l0-NP_hard1}
	B. K. Natarajan, Sparse approximate solutions to linear systems, SIAM journal 
	on computing, 24(2), 227--234, (1995).
	
	\bibitem{l2-NMF}
	V.P. Pauca, J. Piper, R.J. Plemmons, non-negative matrix 
	factorization for spectral data analysis, Linear algebra and its applications, 416(1), 29--47, (2006).
	
	
	
	
	\bibitem{mvnt}
	Y. Y. Tang, Y. Qian, F. Xiong, S. Zeng, J. Zhou, Matrix-Vector 
	non-negative Tensor Factorization for Blind Unmixing of Hyperspectral Imagery, IEEE Transactions on Geoscience and Remote Sensing, 55(3), 1776-1792, (2016).




















\end{thebibliography}
\end{document}